\documentclass[a4paper,12pt,reqno]{amsart}

\usepackage[letterpaper, hmargin=1in, top=1in, bottom=1.2in, footskip=0.7in]{geometry}

\allowdisplaybreaks

\usepackage{rotating,pdflscape,color,amssymb,subfigure,psfrag,amsmath,eufrak,bbm,epsfig}
\usepackage{graphicx}
\usepackage{a4wide,amscd}
\usepackage{tikz} 
 \usepackage[usenames,dvipsnames]{pstricks}
 \usepackage{pst-grad} 

\definecolor{vdarkred}{rgb}{0.6,0,0.2}
\definecolor{vdarkblue}{rgb}{0,0.2,0.6}
\usepackage[pdftex, colorlinks, linkcolor=vdarkblue,citecolor=vdarkred]{hyperref}
\usepackage[small]{caption}

\graphicspath{{./Figures/}}

 \newcommand{\bV}{\mathbf{V}}
\newcommand{\bU}{\mathbf{U}}
 
\newcommand{\bT}{\mathbf{T}}
\newcommand{\bD}{\mathbf{D}}
 
\newcommand{\G}{\mc{G}}

\newcommand{\Cc}{\mc{C}}

\newcommand{\sub}{\substack}

\newcommand{\la}{\label}
 \newcommand{\eqre}{\eqref}
\newcommand{\re}{\ref}
\newcommand{\ld}{\ldots}
\newcommand{\beg}{\begin}
\newcommand{\en}{\end}

\newcommand{\trm}{\textrm}
 
\newcommand{\bgt}{\begin{itemize}}
\newcommand{\ent}{\end{itemize}}
\newcommand{\ite}{\item}

\newcommand{\op}{\operatorname}

\newcommand{\rfl}{\rfloor}
\newcommand{\lfl}{\lfloor}
\newcommand{\si}{\sigma}

\newcommand{\diag}{\operatorname{diag}}

\newcommand{\cat}{\operatorname{Cat}}

\newcommand{\ds}{\displaystyle}

\newcommand{\p}{\mathbb{P}}

\newcommand{\tr}{ \operatorname{Tr}}
\newcommand{\Tr}{\operatorname{Tr}}

\newcommand{\Ninf}{\underset{N\to\infty}{\longrightarrow}}

\newcommand{\Ec}[1]{\E \!\left[ #1 \right]}

\newcommand{\parmo}[4]{\partial_X^{#1}\partial_{\ol X}^{#2}\partial_Y^{#3}\partial_{\ol Y}^{#4} \Phi (\bt)\big|_{\bt=0}}

\newcommand{\OO}[1]{O\! \left( #1 \right)}
\newcommand{\oo}[1]{o\!\left( #1 \right)} 
\newcommand{\mtsmall}{\fontsize{10pt}{8pt}\selectfont}
\newcommand{\vcup}{\mbox{\mtsmall{$\bigcup$}}}

\newcommand{\E}{\mathbb{E}\,}

\newcommand{\R}{\mathbb{R}}
\newcommand{\C}{\mathbb{C}}

\newcommand{\N}{\mathbb{N}}
\newcommand{\z}{\mathbb{Z}}

\newcommand{\ud}{\mathrm{d}}

\newcommand{\wg}{\operatorname{Wg}}
\newcommand{\Moeb}{\operatorname{Moeb}}
\newcommand{\card}{\operatorname{Card}}

\DeclareMathOperator{\eloi}{\overset{\operatorname{law}}{=}}

\DeclareMathOperator{\cloi}{\overset{(d)}{\longrightarrow}}

\DeclareMathOperator{\EE}{\mathcal{E}}
\DeclareMathOperator{\tto}{\longrightarrow}
\newcommand{\Idis}{\operatorname{I}^{\neq}}

\newcommand{\pro}{probability }

\newcommand{\f}{\frac}
\newcommand{\ff}{\frac{1}}
\newcommand{\lf}{\left}
\newcommand{\ri}{\right}
 
\newcommand{\st}{such that }

\newcommand{\lam}{\lambda}

\newcommand{\ti}{\times}

\newcommand{\vfi}{\varphi}
\newcommand{\ste}{\, ;\, }
\newcommand{\mc}{\mathcal }

\newcommand{\eps}{\varepsilon}

\newcommand{\NN}{\mc{N}}

\newcommand{\al}{\alpha}
\newcommand{\tta}{\theta}

\newcommand{\eqlaw}{\stackrel{\textrm{law}}{=}}

\newcommand{\ovl}{\overline}
\newcommand{\ol}{\overline}

\newcommand{\bbm}{\begin{bmatrix}}
\newcommand{\ebm}{\end{bmatrix}}
\newcommand{\bes}{\begin{equation*}}
\newcommand{\ees}{\end{equation*}}
\newcommand{\be}{\begin{equation}}
\newcommand{\ee}{\end{equation}}
\newcommand{\beqy}{\begin{eqnarray}}
\newcommand{\eeqy}{\end{eqnarray}}
\newcommand{\beq}{\begin{eqnarray*}}
\newcommand{\eeq}{\end{eqnarray*}}
\newcommand{\one}{\mathbbm{1}}
\newcommand{\lto}{\longrightarrow}

\newcommand{\ie}{\emph{i.e. }}

\newcommand{\bpm}{\begin{pmatrix}}
\newcommand{\epm}{\end{pmatrix}}

\newcommand{\cd}{\cdots}

\newcommand{\wt}{\widetilde}

\newcommand{\bpr}{\beg{proof}}
\newcommand{\epr}{\en{proof}}

\newcommand{\bet}{\beta}
\newcommand{\del}{\delta}

\newcommand{\pa}{\partial}

\newcommand{\opa}{\op{a}}

\newcommand{\opb}{\op{b}}

\newcommand{\bfa}{\mathbf{a}}
\newcommand{\bfb}{\mathbf{b}}

\newcommand{\bX}{\mathbf{X}}
\newcommand{\bA}{\mathbf{A}}
\newcommand{\bB}{\mathbf{B}}
\newcommand{\bC}{\mathbf{C}}

\newcommand{\bN}{\mathbf{N}}
\newcommand{\bK}{\mathbf{K}}

\newcommand{\bW}{\mathbf{W}}
\newcommand{\bP}{\mathbf{P}}
\newcommand{\bt}{\mathbf{t}}

\newcommand{\bM}{\mathbf{M}}
\newcommand{\bMM}{\mathbf{M}}
\newcommand{\bI}{\mathbf{I}}
\newcommand{\bi}{\mathbf{i}}
\newcommand{\bj}{\mathbf{j}}
\newcommand{\bk}{\mathbf{k}}
\newcommand{\bl}{\mathbf{l}}

 \newcommand{\ii}{\mathrm{i}}

\newcommand{\ka}{\kappa}

  \newcommand{\mre}{\mathrm{e}}
    
    \newcommand{\Circ}{\op{Circle}}
     
    \newcommand{\oti}{\otimes}
 
\newcommand{\mset}[1]{\left\{#1 \right\}_m }
\newcommand{\rset}[1]{\left\{#1 \right\} }

\DeclareMathOperator{\real}{Re}
\DeclareMathOperator{\im}{Im}

\newtheorem{Th}{Theorem}[section]

\newtheorem{propo}[Th]{Proposition}

\newtheorem{lemme}[Th]{Lemma}
\newtheorem{lem}[Th]{Lemma}

\newtheorem{cor}[Th]{Corollary}
\theoremstyle{definition}
 
\newtheorem{rmk}[Th]{Remark}

\newtheorem{Def}[Th]{Definition}

\long\def\symbolfootnote[#1]#2{\begingroup
\def\thefootnote{\fnsymbol{footnote}}\footnote[#1]{#2}\endgroup}

\title{Fluctuations for analytic test functions in the Single Ring Theorem}
\author{Florent Benaych-Georges and Jean Rochet}

\thanks{FBG and JR: MAP5,
Universit\'e Paris Descartes,
45, rue des Saints-P\`eres
75270 Paris Cedex 06, France. florent.benaych-georges@parisdescartes.fr, jean.rochet@parisdescartes.fr.}



\begin{document}

\maketitle

\begin{abstract}
We consider a   non-Hermitian random matrix $\bA$ whose distribution is invariant under the left and right actions of the unitary group. The so-called  \emph{Single Ring Theorem}, proved by Guionnet, Krishnapur and Zeitouni \cite{GUI}, states that the empirical eigenvalue distribution of $\bA$ converges to a limit measure supported by an annulus $S$. In this text, we establish the convergence in distribution of random variables of the type $\tr (f(\bA)\bM)$ where $f$ is analytic on $S$ and the Frobenius norm of $\bM$ has order $\sqrt{N}$. As corollaries, we obtain   central limit theorems for linear spectral statistics of $\bA$ (for analytic test functions) and for finite rank projections of $f(\bA)$ (like matrix entries). As an application, we locate outliers in multiplicative perturbations of $\bA$.
\end{abstract}



\section{Introduction}
The Single Ring Theorem, by Guionnet, Krishnapur and Zeitouni \cite{GUI}, describes the   empirical distribution of the eigenvalues of a large generic matrix with prescribed singular values, \ie an $N\ti N$  matrix of the form $\bA=\bU\bT\bV$, with $\bU, \bV$ some independent Haar-distributed unitary matrices and $\bT$ a deterministic matrix whose singular values are the ones prescribed.    More precisely, under some technical hypotheses, 
 as the dimension $N$ tends to infinity, if  the empirical distribution of the  singular values of $\bA$ converges to a compactly supported limit measure $\Theta$ on the real line, 
then the empirical eigenvalues distribution of $\bA$ converges to a limit measure $\mu$ on the complex plane
which depends only on $\Theta$. The limit measure $\mu$ (see Figure \ref{Fig11})
is rotationally invariant in $\C$ and its support   is the annulus $S:=\{z\in \C\ste a\le |z|\le b\}$, with $a,b\ge 0$ \st \be\la{9914}a^{-2}=\int x^{-2}\ud\Theta( x)\qquad \trm{ and }\qquad b^{2}=\int x^{2}\ud\Theta(  x).\ee  
In this text, we   consider such a matrix $\bA$ 
and we study (Theorem \re{theoremgen}) the joint weak convergence,  as $N\to\infty$, of random variables of the type $$\Tr (f(\bA)\bM),$$ for $f$ an analytic function on the annulus $S$ whose Laurent series expansion has null constant term and $\bM$ a deterministic $N\ti N$ matrix satisfying some limit conditions. These limit conditions (see \eqre{3031512h}) allow to consider both:
 \bgt\ite[--] fluctuations, around their limits as predicted by the Single Ring Theorem, of linear spectral  statistics of $\bA$ (for $\bM=\bI$): $$\Tr f(\bA)=\sum_{i=1}^Nf(\lam_i),$$ where $\lam_1, \ld, \lam_N$ denote the eigenvalues of $\bA$,\\ \ite[--]  finite rank projections  of $f(\bA)$ (for $ \bM=\sqrt{N}\ti$(a matrix with bounded rank)), like the matrix entries of $f(\bA)$. \ent
Let us present  both of these directions with more details.

\subsection{Linear spectral statistics of $\bA$}
As far as Hermitian random matrices are concerned,     linear spectral statistics fluctuations usually come right after the  macroscopic behavior, with the microscopic one, in the natural questions that arise   (see e.g., among the wide literature on the subject,    \cite{jonsson,KKP96,sinai,johansson,BaiYaoBernoulli2005,BAI2009EJP,lytova,MShcherbina11,greg-ofer,lytova,baiysilver,chatterjee,AFCTCL}). For unitary or orthogonal matrices, also, 
many results have been proved (see e.g. the results of Diaconis \emph{et al} in  \cite{dia-shah,dia-evans}, the ones of Soshnikov in \cite{soshni00} or the ones of L\'evy and  Ma\"{\i}da in \cite{thierrymylene}).
For non-Hermitian matrices,   established results are way less numerous: the first one   was \cite{RiderJack}, by Rider and Silverstein,  for analytic test functions  of matrices with i.i.d. entries, then came the paper \cite{RiderVirag} by Rider and Vir\'ag, who managed, thanks to the explicit determinantal structure of the correlation functions of the Ginibre ensemble, to study the fluctuations of linear spectral statistics of such matrices for $\Cc^1$ test functions. Recently, in \cite{ORR}, O'Rourke and  Renfrew studied the fluctuations of linear spectral statistics of elliptic matrices for  analytic test functions and, in \cite{Kemp,CebronKemp}, C\'ebron and Kemp used a dynamical approach to study such fluctations on $\mathbb{GL}_N$. The reason why, except for the breakthrough of Rider and Vir\'ag in \cite{RiderVirag}, many results are limited to analytic test functions is that  
 when non-normal matrices are concerned, functional calculus makes sense only for  analytic functions: if one denotes by $\lam_1, \ld, \lam_N$ the eigenvalues of a non-Hermitian   matrix $\bA$, one can estimate $\sum_{i=1}^N f(\lam_i)$ out of  the numbers $\Tr \bA^k$ or $\Tr((z-\bA)^{-1})$ only when $f$ is analytic. For a $\Cc^2$ test function $f$, one relies on the explicit joint distribution of the $\lam_i$'s or on   Girko's     so-called \emph{Hermitization technique}, which expresses the empirical spectral measure of    $\bA$ as the Laplacian of the function $ z\longmapsto \log|\det(z-\bA)|$ (see e.g. \cite{girko,BOR1}). This is a way more difficult problem, which we consider   in a forthcoming project.
 
 In this text, as a corollary of our main theorem,  we   prove   that for $\bA=\bU\bT\bV$ an $N\ti N$ matrix of the type introduced above and $f$ an analytic function on a neighborhood of the limit support $S$ of the empirical eigenvalue distribution of $\bA$, 
the random variable $$\Tr f(\bA)-\E f(\bA)$$
converges in distribution, as $N\to\infty$, to a centered  complex Gaussian random variable with  a given covariance matrix (see Corollary \re{cor1}). This is a first step in the study of the noise in the Single Ring Theorem. We notice a quite common fact in random matrix theory:  the random variable \bes\la{2903171}
\Tr f(\bA)-\E \Tr f(\bA)=\sum_{i=1}^N f(\lam_i)-\E f(\lam_i)\ees does not need to be renormalized to have a limit in distribution, which reflects the eigenvalue repulsion phenomenon (indeed, would the $\lam_i$'s have been i.i.d., this random variable would have had order $\sqrt{N}$).

Next, two corollaries are given (Corollaries \re{bergman} and \re{carpol174151}), one about the Bergman kernel and the resolvant and one about the log-correlated  limit distribution of the characteristic polynomial out of the support.   
 
 It should be noted that the class of test functions studied -- $f$ analytic on a neighborhood of the  annulus where the eigenvalues $\lam_i$ locate asymptotically  --  is not rich enough to fully characterize the fluctuations the spectrum. For example, not all smooth functions on the annulus can be approximated by analytic functions. 
 Thus while these results do give insight into the fluctuations, the full study of the fluctuations would have to go beyond the realm of analytic test functions.

 \subsection{Finite rank projections and matrix entries}  A century ago, in 1906, \'Emile Borel proved  in  \cite{b1906}  that, for 
a uniformly distributed point $ (X_1 , \ldots , X_N)$ on the unit Euclidian sphere $\mathbb{S}^{N-1}$, the scaled first coordinate $\sqrt{N}X_1$ converges weakly to the   Gaussian distribution as the dimension $N$ tends to infinity. As  explained in the introduction of the paper \cite{diaconis2003}  of  Diaconis  {\it et al.},    this means that the features of the ``microcanonical" ensemble in a certain model for statistical mechanics (uniform measure on the sphere) are captured by the ``canonical" ensemble (Gaussian measure). Since then, a long list of further-reaching results about the asymptotic normality of  entries of   random orthogonal or unitary  matrices  have been obtained (see e.g.   \cite{diaconis2003,meckes08, meckessourav08,collins-stolz08,jiang06,BENCLT,FloGuiJean}).  

 In this text, as a corollary of our main theorem,  we   prove   that for $\bA=\bU\bT\bV$ an $N\ti N$ matrix of the type introduced above,  $f(z)=\sum_{n\in \z}a_nz^n$  an analytic function on a neighborhood of the limit support $S$ of the empirical eigenvalue distribution of $\bA$ and $\bfa,\bfb$ some unit column vectors, 
the random variables of the type  $$\sqrt{N}(\bfb^*f(\bA)\bfa- a_0\bfb^*\bfa),$$ which can be seen as particular cases of the variables of the type $\Tr (f(\bA)\bM)-\E\Tr (f(\bA)\bM)$ for $M=\sqrt N\bfa\bfb^*$, 
converge  in joint distribution, as $N\to\infty$, to   centered  complex Gaussian random variables with  a given covariance matrix (see Corollary \re{cor3}).   This allows for example to consider matrix entries  of $f(\bA)$, in the vein of the works of  Soshnikov \emph{et al.} for Wigner matrices in  \cite{ORRS,PRS12} (see Corollary \re{cor4} and Remark \re{rmk16415abc}). It also applies to the study of   finite rank perturbations of $\bA$ of \emph{multiplicative type}: the BBP phase transition (named after the authors of the seminal paper \cite{BAI}) is well understood for additive or multiplicative  perturbations ($\wt{\bA}=\bA+\bP$ or $\wt{\bA}=\bA(\bI+\bP)$) of general Hermitian   models (see \cite{SP06,CDFF,BEN1} or \cite{BAI,BEN4}), for additive perturbations of various   non-Hermitian models (see \cite{TAO1,FloJean,ROR12,BORCAP1}), but  multiplicative perturbations of non-Hermitian models were so  far unexplored. In Remark \re{rmk16415aha} and Figure \re{Fig:outliersmult}, we explain briefly how our results allow to enlighten a BBP transition for such perturbations. 
 
 \subsection{Organisation of the paper and proofs}

   In the next section, we state our main theorem (Theorem \re{theoremgen}) and its corollaries. The rest of the paper is devoted to the proof of Theorem \re{theoremgen}, to the proof of Corollary \re{carpol174151} and  to the appendix.
   
 Theorem \re{theoremgen} is proved in three steps. First, we do  a cut-off approximation to replace the analytic functions $f$  in the random variables  $\Tr (f(\bA)\bM)$ by   polynomials.  The estimation of the error term in this cut-off is far from obvious relies on non-asymptotic estimates   from   \cite{BEN} and \cite{FloJean}. Then, to prove the convergence to Gaussian random variables, we perform a moment calculation, using Weingarten calculus (for the asymptotic fine moments of Haar-distributed unitary matrices). Weingarten calculus is the theory, due essentially to Collins and \'Sniady,  of the    joint moments of entries of Haar-distributed unitary matrices. We   summarize the necessary ingredients of the theory in Appendix A. The third step is  the computation of the  limit covariance. 
 
\section{Main result}
Let $\bA$ be a random $N\ti N$ matrix implicitly depending on $N$ \st $\bA=\bU\bT\bV$, with $\bU,\bV,\bT$ independent and   $\bU,\bV$   Haar-distributed on the unitary group. We make the following hypotheses on $\bT$:
\beg{assum}\la{assum:1} As $N\to\infty$,  the sequence $\lf(N^{-1}\Tr \bT\bT^*\ri)^{1/2}$ converges in \pro to a deterministic limit $b>0$ and there is $M<\infty $ \st with \pro tending to one, $\|\bT\|_{\op{op}}\le M$. 
\end{assum}

\beg{assum}\la{assum:2} With the convention $1/\infty=0$ and $1/0=\infty$,  the sequence $$\lf(N^{-1}\Tr ((\bT\bT^*)^{-1})\ri)^{-1/2}$$ converges in \pro to a deterministic limit $a\ge 0$. If $a>0$, we also suppose that there is $M'<\infty$ \st with \pro tending to one,
  $\|\bT^{-1}\|_{\op{op}}\le M'$. 
\end{assum}

The following, seemingly purely technical, assumption, which could possibly be relaxed following Basak and  Dembo's approach of \cite{BD13},  is  made to control tails of Laurent series but can be removed  if the $f_j$'s have finite Laurent expansion, like in Corollary \re{cor2} or in Remark \re{rmk16415abc}. Precisely, we need it to cite some estimates from \cite{GUI2}, where they were proved under this assumption.
\beg{assum}\la{assum:3} There exist a constant $\kappa>0$ such that
\begin{eqnarray*} 
\im(z) \ > \ n^{-\kappa} \implies  N^{-1} \big|\im \Tr( (z-\sqrt{\bT\bT^*})^{-1}) \big| \ \leq \ \ff{\kappa}.
\end{eqnarray*}  
\end{assum}

For $f$ an analytic function on a neighborhood of the annulus $$S:=\{z\in \C\ste a\le |z|\le b\} ,$$ the matrix $f(\bA)$ is well defined with \pro tending to one as $N\to\infty$, as it was proved in \cite{GUI2,BEN}   that the spectrum of $\bA$ is contained in any neighborhood of $S$ with \pro tending to one. We denote the Laurent series expansion, on $S$,  of any such function $f$ by $$f(z)=\sum_{n\in \z}a_n(f)z^n.$$

\begin{Th}\label{theoremgen} 
For each $N \geq 1$, let $\bM_1,\ld,\bM_k $ be $N\ti N$ deterministic matrices such that for all $i,j$, as $N \to \infty$,
\beqy\la{3031512h}
\ff N \tr \bM_i \ \tto \ \tau_i, &\ds \ff N \tr \bM_i \bM_j^*   \ \tto \ \al_{ij}, & \ff N \tr  \bM_i\bM_j   \ \tto \ \bet_{ij}
\eeqy
Let $f_1,\ldots,f_k$ be analytic on a neighborhood of  $S$. Then, as $N \to \infty$, the random vector
\be\la{273151}
\Big( \tr  f_j ( \bA)\bM_j  - a_0( f_j )\tr \bM_j  \Big)_{j=1}^k
\ee
converges to a centered complex Gaussian vector $(\mathcal{G}(f_1),\ldots,\mathcal{G}(f_k))$  whose distribution is defined by 
\beq
\E \mathcal{G}(f_i) \mathcal{G}(f_j) & = &  \sum_{n \geq 1}\big((n-1)\tau_i\tau_j+\bet_{ij} \big)\big( a_n(f_i) a_{-n}(f_j) + a_{-n}(f_i) a_{n}(f_j) \big)\\
\E \mathcal{G}(f_i) \ol{\mathcal{G}(f_j)} & = & \sum_{n \geq 1}\big((n-1)\tau_i\ol{\tau_j}+\al_{ij}\big)\big(a_n(f_i)\ol{a_n(f_j)}b^{2n} + a_{-n}(f_i)\ol{a_{-n}(f_j) } 	a^{-2n}  \big).
\eeq
\end{Th}
\begin{rmk} In \eqre{273151}, $$\tr  f_j ( \bA)\bM_j  - a_0( f_j )\tr \bM_j$$ rewrites $$\tr  f_j ( \bA)\bM_j  - \E\tr  f_j ( \bA)\bM_j.$$
Indeed,  $  \E f(\bA)=a_0 \bI$, as a consequence of the fact that for any $n\ne 0$, for any $\tta\in \R$, $  \bA^n\eqlaw \mre^{\ii\tta}\bA^n$, which follows from   the invariance of the Haar measure.
  \end{rmk}

\begin{rmk}\la{rmk:series_sommables}Note that   if $a=0$, as the $f_j$'s are analytic on $S$, we have   $a_{-n}(f_j)=0$ for each $n\ge 1$ and each $j$, so that the above expression still makes sense. Besides, 
it seems reasonable to verify that the two series above converge:
\beq
\sum_{n \geq 1}n |a_n(f_i)| |a_n(f_j)| b^{2n} & \leq & \big(\max_{n \geq 1} |a_n(f_j)|b^n\big) \sum_{n \geq 1} n |a_n(f_i)|b^n \ < \ \infty  \\
\sum_{n \geq 1}n |a_n(f_i)|| a_{-n}(f_j)| & \leq &\big( \max_{n \geq 1}  |a_n(f_i)|b^n \big) \sum_{n \geq 1} n | a_{-n}(f_j)| a^{-n} \ < \ \infty.
\eeq
\end{rmk}

\beg{rmk}[Relation to second order freeness] A theory has been developed recently about Gaussian  fluctuations (called {\it second order limits}) of traces of large random matrices, the most emblematic articles in this theory   being \cite{mingo-nica04, mingo-speicher06, mingo-piotr-speicher07, mingo-piotr-collins-speicher07}. Theorem \re{theoremgen} can be compared to some of these  results. However, our hypotheses on the matrices we consider are   of a different nature  than the ones of the previously cited papers, since the convergence of the non commutative distributions is not required here: our hypotheses are satisfied    for example by matrices like $\bM_j=\sqrt{N}\times$(an elementary $N\times N$ matrix), which have no bounded moments of order higher than two.
\en{rmk}

Our two main applications are the case where the $\bM_j$'s are all $\bI$ (Corollaries \re{cor1} and \re{cor2}) and the cases where the $\bM_j$'s are   $\sqrt{N}$ times   matrices with bounded rank and norm, like elementary matrices  (Corollaries \re{cor3} and \re{cor4}).
In the case  $\bM = \bI$,   we immediately obtain the following corollary  about   linear spectral statistics of $\bA$.

\begin{cor}\label{cor1}
Let $f_1,\ldots,f_k$ be analytic on a neighborhood of $S$. Then, as $N \to \infty$, the random vector
$$
\Big( \tr f_j \big( \bA\big) - Na_0( f_j )   \Big)_{j=1}^k
$$ 
converges to a centered complex Gaussian vector $(\mathcal{G}(f_1),\ldots,\mathcal{G}(f_k))$ such that
\beq
\E \mathcal{G}(f_i) \mathcal{G}(f_j)  & = &  \sum_{n \in \z}|n| a_n(f_i) a_{-n}(f_j)  \\
\E \mathcal{G}(f_i) \ol{\mathcal{G}(f_j)}  & = & \sum_{n \geq 1}n\big(a_n(f_i)\ol{a_n(f_j)}b^{2n} + a_{-n}(f_i)\ol{a_{-n}(f_j) } 	a^{-2n}  \big).
\eeq
\end{cor}

For $n\ge 1$, let us define the functions $$  \varphi^\pm_n(z) \ := \  \lf(\f{z}{b}\ri)^n \pm   \lf(\f{a}{z}\ri)^n.$$
These functions (plus the constant one)  define a basis of the space of analytic functions on a neighborhood of $S$  and we have the change of basis formula
$$\sum_{n\in \z} a_n z^n=a_0+\sum_{n\ge 1} c_n^+\vfi_n^+(z)+c_n^-\vfi_n^-(z)\iff \forall n\ge 1,\; \bpm a_n\\ a_{-n}\epm = \bpm b^{-n}& b^{-n}\\ a^n& -a^n\epm \bpm c_n^+\\ c_n^- \epm\,,$$ 
 implying that\beq
\sum_{n\geq 1}|a_n(f)|^2b^{2n} + |a_{-n}(f)|^2 	a^{-2n}  & = & 2\sum_{n \geq 1}|c^+_n(f)|^2 + |c^-_n(f)|^2.
\eeq
Besides, these functions allow to identify the underlying white noise in Theorem \re{theoremgen} (we only state it here in the case  $\bM_j=\bI$, but this of course extends to the case of general $\bM_j$'s, allowing for example to state analogous results for the matrix entries).

\beg{cor}[Underlying white noise]\label{cor2}The finite dimensional marginal distributions of  $$(\Tr \vfi^+_n(\bA))_{n\ge 1}\,\vcup \,(\Tr \vfi^-_n(\bA))_{n\ge 1}$$ converge to the ones of  a collection   $(\mathcal{G}_n^+)_{n\ge 1}\, \cup \,(\mathcal{G}_n^-)_{n\ge 1}$ of independent centered complex Gaussian random variables satisfying $$\E |\mathcal{G}_n^\pm|^2\;=\;2\qquad;\qquad \E (\mathcal{G}_n^\pm)^2\;=\;\pm 2n(a/b)^n.$$
\en{cor}
 
\begin{rmk}[Ginibre matrices]
In the particular case where $\bA$ is a Ginibre matrix (\ie with i.i.d. entries with law  $ \NN_\C(0,N^{-1})$), we reproduce  the result of Rider and Silverstein \cite{RiderJack}, noticing that in this case $a=0$ and $b=1$, so that $a_n(f)=0$ when $n<0$ and $\E \mc{G}(f_i)\mc{G}(f_j) =0$, and, for $\ud m(z)$   the Lebesgue measure on $\C$,
\beq
&&\ff{\pi}\int_{|z|<1} \f{\pa}{\pa z}f_i(z) \ol{\f{\pa}{\pa z}f_j(z)}\ud m(z) \\
& = & \ff \pi \int_{|z|<1}-\ff{4\pi^2} \oint_{\Circ(1+\eps)} \oint_{\Circ(1+\eps)} \f{f_i(\xi_1)}{(\xi_1 - z)^2} \f{\ol{f_j(\xi_2)}}{(\ol{\xi_2} - \ol z)^2}\ud \xi_1\ud \xi_2\ud m(z)\\
& = & -\ff{4\pi^2} \oint_{\Circ(1+\eps)} \oint_{\Circ(1+\eps)} \!\!\!\!\f{f_i(\xi_1)\ol{f_j(\xi_2)}}{\xi_1^2 \ol{\xi_2}^2} \ff \pi \int_{|z|<1} \sum_{n,n' \geq 1}nn' \lf(\f{z}{\xi_1}\ri)^{n-1} \lf(\f{\ol z}{\ol \xi_2}\ri)^{n'-1}\!\!\!\!\ud m(z)\ud \xi_1\ud \xi_2\\
& = & -\ff{4\pi^2} \oint_{\Circ(1+\eps)} \oint_{\Circ(1+\eps)} f_i(\xi_1)\ol{f_j(\xi_2)}\sum_{n \geq 1}n \big(\xi_1\ol{\xi_2}\big)^{-n-1}\ud \xi_1\ud \xi_2 \\
& = & \sum_{n \geq 1} n a_n(f_i) \ol{a_n (f_j)} 
\eeq 
\end{rmk}
\begin{rmk}
If $\bT = \bI$, and the $f_k$'s are polynomial, we reproduce  a result of Diaconis and Shahshahani  \cite[Theorem 1]{dia-shah}   on the limit joint distribution of
$$
\lf(\tr(\bU^k)\ri)_{k=1}^n,
$$
 where $\bU$ is Haar-distributed. Actually, the Corollary \ref{cor1} is slightly stronger, since the result holds for $\bA = \bU\bT$ or $\bA = \bU\bT\bV$ as long as $\bT$ satisfies   
\beqy \la{1827020420051}
 \lim_{N \to \infty} \ff N \tr(\bT\bT^*) \ = \lim_{N \to \infty}  \ff N \tr((\bT\bT^*)^{-1})  = 1,
\eeqy
in which case  $\bA$ may be seen as a multiplicative perturbation of $\bU$. Indeed, \eqref{1827020420051} implies that all singular values of $\bT$ are close to $1$. The matrix $\bT$ satisfies the condition \eqref{1827020420051}  for example if it is diagonal and all its diagonal coefficients are equal to $1$ except   $\oo{N}$ of them (which stay away from $0$ and $\infty$).
\end{rmk}

 \beg{cor}[Bergman kernel and resolvant]\la{bergman}The random process $$ \lf(\Tr (z-\bA)^{-1}\ri)_{|z|<a}\,\cup\,\lf(\Tr (z-\bA)^{-1}\ri)_{|z|>b}$$ converges, for the finite-dimensional   distributions, to  a centered complex  Gaussian process  $$\lf(\mc{G}_z\ri)_{|z|<a}\,\cup\, \lf(\mc{H}_z\ri)_{|z|>b}$$ with covariance defined by $$ \E \mc{G}_z\ovl{\mc{G}_{z'}}=\f{b^2}{(b^2-z\ovl{z'})^2}\,,\qquad    \E \mc{H}_z\ovl{\mc{H}_{z'}}=\f{a^2}{(a^2-z\ovl{z'})^2}\,,\qquad \E \mc{G}_z \mc{H}_{z'} =-\ff{(z'-z)^2}.$$ and by the fact that  $$\forall \tta\in \R,\qquad \lf(\mre^{-\ii\tta}\mc{G}_z\ri)_{|z|<a}\,\cup\, \lf(\mre^{\ii\tta}\mc{H}_z\ri)_{|z|>b}\,\,\eqlaw\,\, \lf(\mc{G}_z\ri)_{|z|<a}\,\cup\, \lf(\mc{H}_z\ri)_{|z|>b}.$$
 \en{cor}
 
 \beg{rmk}The   kernel of the limit Gaussian analytic function, in the previous corollary, is, up to a constant factor, the \emph{Bergman kernel} (see \cite{bell,peresvirag}).
 \en{rmk}

 \beg{cor}[Characteristic polynomial out of the support]\la{carpol174151}The random process $$ \lf(\log|\det(z-\bA)|- \Tr\log\bT\ri)_{|z|<a}\,\cup\,\lf(\log|\det(z-\bA)|-N\log|z|\ri)_{|z|>b}$$ converges, for the finite-dimensional   distributions, to  a centered real  Gaussian process  $$\lf(\mc{G}_z\ri)_{|z|<a}\,\cup\, \lf(\mc{H}_z\ri)_{|z|>b}$$ with covariance defined by  $$2\E \mc{G}_z\mc{G}_{z'}=-\log\lf|1-\f{z\ovl{z'}}{a^2}\ri|\,,\qquad 2\E \mc{H}_z\mc{H}_{z'}=-\log\lf|1-\f{b^2}{z\ovl{z'}}\ri|\,,\qquad 2\E \mc{G}_{z}\mc{H}_{z'}=-\log\lf|1-\f{z}{z'}\ri|\,.$$ 
 \en{cor}
 
 \beg{rmk}As $z\ne z'$ both tend to the same point on the boundary of $S$, the above covariances are equivalent to $-\log|z-z'|$.  In the light of the log-correlation approach to   the Gaussian Free Field (see \cite{BRSVoverview}), it supports the idea that  on the limit support $S$, the characteristic polynomial of $\bA$ should tend to an object related to the  Gaussian Free Field, as for Ginibre matrices (see Corollary 2 of \cite{RiderVirag}). It would be interesting to see to what extent such a convergence   depends on the hypotheses made on the precise distribution of the singular values of $\bT$. 
 \en{rmk}

In the case  $\bM_j = \sqrt{N}\bfa_j\bfb^*_j$, we immediately obtain the following corollary:
\begin{cor}\label{cor3}For each $N \geq 1$, let $\bfa_1,\bfb_1,\ld,\bfa_k ,\bfb_k$ be  deterministic  column vectors with size $N$   such that for all $i,j$, as $N\to\infty$, 
\be\la{3031512huv}
  \bfa_i^*\bfa_j \ \tto \ \ka^{\opa,\opa}_{ij}\in \C\qquad;\qquad \bfb_i^*\bfa_j \ \tto \ \ka^{\opb,\opa}_{ij}\in \C \qquad;\qquad \bfb_i^*\bfb_j \ \tto \ \ka^{\opb,\opb}_{ij} \in \C
\ee
Let $f_1,\ldots,f_k$ be analytic on a neighborhood of $S $. Then, as $N \to \infty$, the random vector
\be\la{273151uv}
\sqrt{N}\Big(   \bfb_j^*f_j ( \bA)\bfa_j  -  \bfb_j^*  \bfa_j a_0( f_j ) \Big)_{j=1}^k
\ee
converges to a centered complex Gaussian vector $(\mathcal{G}(f_1),\ldots,\mathcal{G}(f_k))$  \st
\beq
\E \mathcal{G}(f_i) \mathcal{G}(f_j)  & = &  \sum_{n \geq 1}  \ka_{ji}^{\opa,\opa}\ka_{ij}^{\opb,\opb}  \big( a_n(f_i) a_{-n}(f_j) + a_{-n}(f_i) a_{n}(f_j) \big)\\
\E \mathcal{G}(f_i) \ol{\mathcal{G}(f_j)} & = & \sum_{n \geq 1} \ka_{ji}^{\opb,\opa}\ka_{ij}^{\opb,\opa}\big(a_n(f_i)\ol{a_n(f_j)}b^{2n} + a_{-n}(f_i)\ol{a_{-n}(f_j) } 	a^{-2n}  \big).
\eeq
\en{cor}

\beg{rmk}[Application to multiplicative finite rank perturbations of $\bA$]\la{rmk16415aha}The previous corollary has several applications to the study of the outliers of spiked models related to the Single Ring Theorem. It allows for example   to understand easily, using the techniques developed in \cite{FloJean}, the impact of \emph{multiplicative} finite rank perturbations on the spectrum of $\bA$ (whereas only \emph{additive} perturbations had been studied so far).
For example, one can deduce from this corollary   that for $\bP$ a deterministic matrix with bounded operator norm and rank one, if one defines $\wt{\bA}:=\bA(\bI+\bP)$ and $\hat{\bA}:=\bA(\bI+\bA\bP)$, then \bgt\ite the matrix $\wt{\bA}$ has no outlier  (\ie the support of its spectrum still converges to $S$), 
\ite  the matrix $\hat{\bA}$ has no outlier  with modulus $>b$, but each eigenvalue $\lam$ of $\bP$ 
    \st $|\lam|>a^{-1}$ gives rise to  an outlier of $\hat{\bA}$ located approximately at $-\lam^{-1}$ (besides, when the multiplicity of $\lam$ as an eigenvalue of $\bP$ is $1$, the fluctuations of the outlier  around $-\lam^{-1}$ are Gaussian and with order $N^{-1/2}$).
\ent
This phenomenon is illustrated by Figure \re{Fig:outliersmult}.
\newcommand{\scalefigoutliersmult}{.45}
\begin{figure}[ht]
\centering
\subfigure[Spectrum of $\bA$]{
\includegraphics[scale=\scalefigoutliersmult]{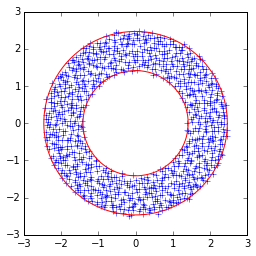} 
\label{Fig11}} \qquad 
\subfigure[Spectrum of $\wt{\bA}:=\bA(\bI+\bP)$]
{\includegraphics[scale=\scalefigoutliersmult]{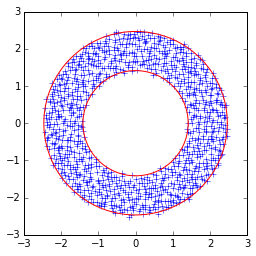} 
\label{Fig12}}
\qquad 
\subfigure[Spectrum of $\hat{\bA}:=\bA(\bI+\bA\bP)$ (small circles are centered at the theoretical limit locations of the outliers)]
{\includegraphics[scale=\scalefigoutliersmult]{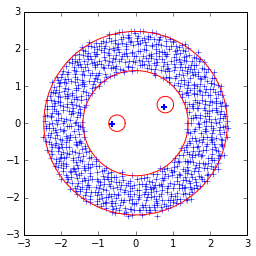} 
\label{Fig13}}
\caption{{\bf Outliers/lack of outliers for multiplicative perturbations:} simulation realized with a single  $10^3\times 10^3$  matrix $\bA=\bU\bT\bV$   when the singular values of $\bT$ are uniformly distributed on $[0 . 5 , 4]$  and $\bP=\diag(-2,(0.8+0.5\ii)^{-1},1/3,0, \ld, 0)$. As predicted, none of these matrices has any outlier outside   the outer circle, nor do the two first ones inside the inner circle, but $\hat{\bA}$ has two outliers  inside the inner circle, close to  the predicted locations.}\la{Fig:outliersmult}
\end{figure}
\end{rmk}

To state the next corollary,   let us first give the definition of Gaussian elliptic matrices.

\begin{Def}[Gaussian elliptic matrices]\la{defelliptic}
Let $\rho\in \C$ \st $|\rho| \leq 1$.   A \emph{Gaussian elliptic matrix with parameter $\rho$} is an $N \ti N$  Gaussian centered complex   random matrix $\bX = (x_{ij})$ satisfying :
\bgt
\ite[(1)] the random vectors  $\big(x_{ij},x_{ji} \big)_{i \le  j}$ are  independent,
\ite[(2)] $\forall i$, $\E|x_{i i}|^2 = 1$ and $\E x_{i i}^2  = \rho$,
\hspace{0.1mm}
\ite[(3)] $\forall i\neq j$, $\E{|x_{i j}|^2}=1$, $\E x_{i j}^2  = 0$, $\E x_{i j}x_{j i} = \rho$ and $\E x_{i j}\ol{x_{j i}}  = 0$. 
\ent
\end{Def}

This matrix ensemble was introduced by Girko in \cite{girko}, and its name is due to the fact that its empirical eigenvalue distribution is the uniform distribution inside an ellipse. In the case where $\rho=0$ (resp. $\rho=1$), we get a Ginibre (resp. GUE) matrix.

\begin{cor}\label{cor4}
Let $f$ be analytic on a neighborhood of $S  $  such that
$$
\sum_{n \geq 1} |a_n(f)|^2 b^{2n} + |a_{-n}(f)|^2 a^{-2n} \ = \ 1.
$$
Let $k$ be a  fixed positive integer  and let $I=I(N)$ be a (possibly $N$-dependent) subset of $\{1, \ld, N\}$ with cardinality $k$.  Let us define the random $k\ti k$ matrix  $$\bX_N:=\sqrt{N}\bbm f(\bA)_{ij}-a_0(f)\del_{ij}\ebm_{(i,j)\in I\ti I}.$$  Then, as $N \to \infty$, the matrix $\bX_N$ converges in distribution to a $k\ti k$ Gaussian elliptic matrix $\bX$ with parameter $\rho$,
 $$
 \rho \ := \  2\sum_{n \geq 1}a_n(f)a_{-n}(f).
 $$ 
\end{cor}

\begin{rmk}\la{rmk16415abc}\bgt\ite[a)]  In the case where $f(z)=z$, 
we rederive the well-known result that any fixed-size principal submatrix  of $\sqrt{N}\bU$   converges to a Ginibre matrix (see e.g. \cite{diaconis2003,jiang06}). 
\ite[b)] By Corollary \re{cor4},  the statement of the first part of this remark happens to stay true, up to a constant multiplicative factor,  if $\bU$ is replaced by $\bA=\bU\bT\bV$ or even by $\bA^n$  or by $f(\bA)$ if $f$ is analytic in a neighborhood of the disc $\ovl{B}(0,b)$.   
\ite[c)] It also follows from what precedes that for any $n\ge 1$, any sequence of principal submatrices with fixed size of $\sqrt{N/2}(\bU^n+\bU^{-n})$ and $\sqrt{N/2}(\bU^n-\bU^{-n})$  converge in distribution to   a GUE matrix and $\ii$ times a GUE matrix, both being independent.\ent
\end{rmk}

\section{Proof of Theorem \re{theoremgen}}

To avoid having to treat the cases $a>0$ and $a=0$ separately all along the proof,    we shall suppose that $a>0$ (the case   $a=0$   is more simple, as sums run only on $n\ge 0$). Besides, note that by invariance of the Haar measure, the distribution of the random matrix $\bA$ depends on $\bT$ only through its singular values, so we shall  suppose that $\bT=\diag(s_1, \ld, s_N)$, with $s_i\ge 0$.
At last, as the limit distributions, in Theorem \re{theoremgen},  only depend on $\bT$ only through the deterministic parameters $a,b$, up to a conditioning, one can suppose that $\bT$ is deterministic (and that both $\|\bT\|_{\op{op}}$ and $\|\bT^{-1}\|_{\op{op}}$ are uniformly bounded, by Asssumptions \re{assum:1} and \re{assum:2}). 

\subsection{Randomization of the $\bM_j$'s}\la{sec:randN}Let us define $\bW:=\bV\bU$. The random matrix $\bW$ is also Haar-distributed and independent from $\bV$. Besides, for each $j$, as $\bA=\bU\bT\bV= \bV^*\bW\bT\bV$, $$\Tr f_j(\bA)\bM_j= \Tr \bV^*f_j(\bW\bT)\bV\bM_j=\Tr f_j(\bW\bT)\bV\bM_j\bV^*$$
As a consequence, we shall suppose that  $\bA=\bU\bT$ (instead of $\bA=\bU\bT\bV$) and   that there is a Haar-distributed random unitary matrix $\bV$, independent of $\bU$, \st for each $j$, $\bM_j=\bV\wt{\bM}_j\bV^*$, with $\wt{\bM}_1, \ld, \wt{\bM}_k$ a collection of deterministic matrices also satisfying \eqre{3031512h}.

\subsection{Tails of the series}\la{section152327042015}  Let us first prove that Theorem \re{theoremgen} can be deduced from the particular case where there is $n_0$ \st for all $n$, we have $$|n|>n_0\implies \forall j=1, \ld, k, \; a_n(f_j)=0.$$

 Let $\eps\in (0,a/2)$ \st the domain of each $f_j$ contains the annulus of complex numbers $z$ \st  $ a-2\eps\le|z|\le b+2\eps$.

\beg{lem}\la{lem:3031514h} 
There is a constant $C$ independent of $N$ \st for any $n$ \st $n^6\le N$ and any $j=1, \ld, k$, we have 
$$ \E  |\Tr \bA^n\bM_j|^2 \le C n^2\lf(\one_{n\ge 0}   (b+\eps)^{2n}+\one_{n\le 0}   (a-\eps)^{2n}\ri)$$
\end{lem}
\bpr
With the notation of Section \re{sec:randN}, let $\E_\bV$ denote the expectation with respect to the randomness of $\bV$.  For each $n\in \z$ and each $j$, by Lemma \re{lem161423031500}, we have 
\beq \E_{\bV} |\Tr \bA^n\bM_j|^2&=&\E_\bV \Tr \bA^n\bV\wt{\bM}_j\bV^*\Tr  (\bA^*)^n\bV\wt{\bM}^*_j\bV^*\\
&=&\ff{N^2-1}\lf(\Tr \bA^n\Tr(\bA^*)^n\Tr\wt{\bM}_j\Tr\wt{\bM}^*_j+ \Tr\bA^n (\bA^*)^n\Tr \wt{\bM}_j\wt{\bM}_j^*\ri)\\
& &-\ff{N(N^2-1)}\lf(\Tr \bA^n\Tr(\bA^*)^n\Tr\wt{\bM}_j \wt{\bM}^*_j+\Tr\bA^* (\bA^*)^n\Tr \wt{\bM}_j\Tr\wt{\bM}_j^*\ri)\\
&\le & \ff{N^2-1}\lf(|\Tr \bA^n|^2\Tr\wt{\bM}_j\Tr\wt{\bM}^*_j+ \Tr\bA^n (\bA^*)^n\Tr \wt{\bM}_j\wt{\bM}_j^*\ri)\\
&\le & C\lf(|\Tr \bA^n|^2 +N^{-1}\Tr\bA^n (\bA^*)^n \ri),
\eeq where $C$ is a constant independent of $N$. Then the conclusion follows from Lemma \re{Th1JAF}.
\epr

\beg{lem}\la{lem:3031514h1}There are some constants $C>0$ and $c\in (0,1)$   and  a sequence $\EE=\EE_N$ of events \st $$\p(\EE )\Ninf 1$$ and \st  for all $N$, all $n_1\ge 0$ and all $j=1, \ld, k$, we have
$$\E\Big|\one_{\EE}\sum_{|n|>n_1}  a_n(f_j) \Tr (\bA^n\bM_j) \Big|\;\le \; CN(1-c)^{n_1}.$$
\end{lem}

\bpr  By \cite[Lem. 3.2]{FloJean}, we known that there is a constant $C_1$ \st  the event $$\EE=\EE_N:=\{\forall n\ge 0 , \, \|\bA^n\|_{\op{op}}\le C_1(b+\eps)^n\}\cap \{\forall n\le 0 , \, \|\bA^{ n}\|_{\op{op}}\le C_1(a-\eps)^{ n}\}$$ has \pro tending to one.

Then one concludes easily, noting first that by non-commutative H\"older inequalities (see \cite[Eq. (A.13)]{agz}), we have $$\one_{\EE}|\Tr (\bA^n\bM_j)|\le \beg{cases} C_1(b+\eps)^nN\sqrt{N^{-1}\Tr \bM_j\bM_j^*}&\trm{ if $n\ge 0$}\\ \\
C_1(a -\eps)^{ n}N\sqrt{N^{-1}\Tr \bM_j\bM_j^*}&\trm{ if $n\le 0$}\en{cases}$$ and secondly that there is $c\in (0,1)$ \st for each $j$,  the sequences $$\lf( a_n(f_j)  \f{(b+\eps)^n}{(1-c)^n}\ri)_{n\ge 0}  \qquad ;\qquad \lf( a_n(f_j)  \f{(a-\eps)^{n}}{(1-c)^{-n}}\ri)_{n\le 0}$$ are bounded.
\epr  

As a consequence of Lemmas \re{lem:3031514h} and \re{lem:3031514h1}, for any $0<n_0<n_1\le N^{1/6}$ and any $j=1, \ld, k$, 
\beq \E\Big|\one_{\EE}\sum_{|n|>n_0}  a_n(f_j) \Tr (\bA^n\bM_j) \Big|&\le & \sum_{n_0<|n|\le n_1}  |a_n(f_j)|\sqrt{\p(\EE)}\sqrt{\E |\Tr (\bA^n\bM_j)|^2}\\ &&+\E\Big|\one_{\EE}\sum_{|n|>n_1}  a_n(f_j) \Tr (\bA^n\bM_j) \Big|
\\
&\le & \sum_{n_0<|n|\le n_1}C |a_n(f_j)| n^2\lf(\one_{n\ge 0}   (b+\eps)^{2n}+\one_{n\le 0}   (a-\eps)^{2n}\ri)\\ &&+CN(1-c)^{n_1} 
\eeq
Choosing first $n_1=\lfl A\log N\rfl$ for $A$ a large enough constant and then using the fact that for any $j=1, \ld, k$, $$\sum_{n\in \z}  |a_n(f_j)| n^2\lf(\one_{n\ge 0}   (b+\eps)^{2n}+\one_{n\le 0}   (a-\eps)^{2n}\ri)<\infty,$$we deduce that for any $\del>0$, there is $n_0>0$ fixed \st for all $N$ large enough, \be\la{eq:tail} \sum_{j=1}^k\E\Big|\one_{\EE}\sum_{|n|>n_0}  a_n(f_j) \Tr (\bA^n\bM_j) \Big|\;\le \; \del,\ee for $\EE=\EE_N$ as in Lemma \re{lem:3031514h1}.
 
Let us now suppose Theorem \re{theoremgen} to be proved in  the particular case where there is $n_0$ \st for all $n$, we have $$|n|>n_0\implies \forall j=1, \ld, k,\; a_n(f_j)=0$$ and let us prove it in the general case. 
Let $X_N$ denote the random vector of  \eqre{273151}. We want to prove that  as $N\to\infty$, the distribution of $X_N$ tends to the one of $\G:=(\mathcal{G}(f_1),\ldots,\mathcal{G}(f_k))$, \ie that for any function $F:\C^k\to \C$ which is $1$ Lipschitz and bounded by $1$, we have $$\E F(X_N)\Ninf \E F(\G).$$
To do so, we first set
\beq
X_{N,n_0} & := & \Big( \sum_{|n|<n_0}a_n(f_j)\tr \big(   \bA^n\bM_j\big) - \tr(\bM_j)a_0(f_j) \Big)_{j=1}^k
\eeq
By hypothesis, for any fixed $n_0$,    $X_{N,n_0}$ converges in distribution  to a  centered complex Gaussian vector $\mathcal{G}_{n_0}:=(\mathcal{G}_{n_0}(f_1),\ldots,\mathcal{G}_{n_0}(f_k))$ such that
\beq
\Ec{\mathcal{G}_{n_0}(f_i) \mathcal{G}_{n_0}(f_j)} & = & \Ec{\mathcal{G}(f_i) \mathcal{G}(f_j)} + \eta^{n_0}_{ij} \\
\Ec{\mathcal{G}_{n_0}(f_i) \ol{\mathcal{G}_{n_0}(f_j)}} & = & \Ec{\mathcal{G}(f_i) \ol{\mathcal{G}(f_j)}} + \delta^{n_0}_{ij},
\eeq
where $\ds\lim_{n_0 \to \infty}\sum_{1 \leq i,j \leq k} |\eta^{n_0}_{ij}|+|\delta_{ij}^{n_0}| = 0$.
Therefore,
\beq
&&\lf|\Ec{F(X_N) - F(\mathcal G)} \ri| \\ & \leq & \lf|\Ec{F(X_N) - F(X_{N,n_0})} \ri| + \lf|\Ec{F(X_{N,n_0}) - F(\mathcal{ G}_{n_0})} \ri| +\lf|\Ec{ F(\mathcal{ G}_{n_0})- F(\mathcal G)} \ri| \\
& \leq &  2 \p(\EE^c) + \Ec{\one_{\EE}\|X_N - X_{N,n_0}\|_2^2} +  \lf|\Ec{F(X_{N,n_0}) - F(\mathcal{ G}_{n_0})} \ri| + \Ec{\|\mathcal{ G}_{n_0}- \mathcal G\|^2_2}
\eeq
which can be as small as we want by \eqref{eq:tail} and  the fact that $X_{N,n_0} \cloi \mathcal{G}_{n_0}$   if $\mathcal{ G}_{n_0}$ and $\mathcal{ G}$ are coupled the right way.

\subsection{Proof of Theorem \ref{theoremgen} when the $f_j$'s are polynomial in $z$ and $z^{-1}$} We suppose here  that there is   $n_0>0$ such that for all $n > n_0$ and all $1 \leq j \leq k$, $a_n(f_j)=0$. Without any loss of generality, we also assume that for all $j$, $a_0(f_j)=0$. In this case, any linear combination of the $ \Tr f_j(\bA)\bM_j $'s can be written 
\beq
 G_N & := & \sum_{j=1}^{k} \nu_{j}\tr f_j(\bA)\bM_j
\; = \;  \sum_{|n|\leq n_0} \tr \bA^n \bN_n
\eeq
where $\ds \bN_n := \sum_{j=1}^{k} \nu_j a_n(f_j) \bM_j$. 
  Written this way, we notice that to prove that the limit distribution of $G_{N}$ is Gaussian, we simply have to prove that the random vector 
$$
\left( \tr \bA^n \bN_n   \right)_{-n_0 \leq n \leq n_0}
$$
converges in distribution  to a   Gaussian vector. We will prove it   by computing the limit of the joint moments.

 Before going any further, recall that 
by the preliminary randomization of the $\bN_j$'s from section \re{sec:randN},  we suppose that $\bA=\bU\bT$ (instead of $\bA=\bU\bT\bV$) and   that there is a Haar-distributed random unitary matrix $\bV$, independent of $\bU$, \st for each $j$, $\bN_j=\bV\wt{\bN}_j\bV^*$, with $\wt{\bN}_j$ a     deterministic matrix.

We shall proceed in three steps:
\bgt
\ite[{\bf a)}]  \ First, we prove the asymptotic independence of the random vectors $$\lf(\Tr \bA^n\bN_n,\Tr \bA^{-n}\bN_{-n}\ri)_{n\ge 1} $$
by the factorization of the joint moments.  More precisely, we prove, thanks to   Corollary \ref{lemsauveur15460502},  that for any $(p_{n})_{n=1}^{n_0}$, $(q_{n})_{n=1}^{n_0}$, $(r_{n})_{n=1}^{n_0}$, $(s_{n})_{n=1}^{n_0}$,
\beq
&&\E\Big[\prod_{1 \leq n \leq n_0} \big(\tr \bA^{n}\bN_n \big)^{p_{n}}\ol{\big(\tr \bA^{n}\bN_n \big)^{q_{n}}} \big(\tr \bA^{-n}\bN_{-n} \big)^{r_{n}}\ol{\big(\tr \bA^{-n}\bN_{-n} \big)^{s_{n}}}\Big] \\
&=&\prod_{1 \leq n \leq n_0}\E\Big[ \big(\tr \bA^{n}\bN_n \big)^{p_{n}}\ol{\big(\tr \bA^{n}\bN_n \big)^{q_{n}}} \big(\tr \bA^{-n}\bN_{-n} \big)^{r_{n}}\ol{\big(\tr \bA^{-n}\bN_{-n} \big)^{s_{n}}}\Big]  +\oo1
\eeq
\ite[{\bf b)}]  \ Then, we prove for any  fixed $n$, the random complex vector $$(\tr \bA^n  \bN_n ,\tr \bA^{-n}  \bN_{-n} )$$ converges in distribution to a  centered complex Gaussian vector thanks to the criterion  provided by the Lemma \ref{150112032015}. This criterion consists in  proving that the joint moments, at  the large $N$ limit,   satisfy the same induction relation as the moments of a Gaussian   distribution.
\ite[{\bf c)}]  \ It will follow from {\bf a)} and {\bf b)} that when all $f_j$'s are polynomials in $z$ and $z^{-1}$, the random vector of \eqre{273151} converges in distribution to a centered Gaussian vector. 
To conclude the proof, the last step will be to prove that the limit  covariance is the one given in Theorem \re{theoremgen}.
\ent

\noindent In the proofs of {\bf a)} and {\bf b)}, we shall need to compute expectations with respect to the randomness of the Haar-distributed matrix $\bU$. More precisely, we shall need to compute sums of expectations with respect of $\bU$ resulting from the expansion   of products of traces involving powers of $\bA$ (such as $\tr \bA^n \bN_n $). To do so, we will use the Weingarten calculus (see Proposition \ref{wg}) and  shall always proceed in the following way: first, we use \eqref{wg1} to state that all the terms of the sum are null except those for which the left (resp. right) indices involved in $u$ are obtained by permuting the left (resp. right) ones involved in $\ol{u}$. Then, we claim, by Remark \ref{alldistinct},  that among the remaining terms, we can neglect all   those whose indices are not pairwise distinct. At last, once all the remaining terms are, up to multiplicative constant, equal to $\wg(\sigma)$ for some permutation $\sigma$, we neglect all those for which $\sigma \neq id$ (see Remark \ref{moebius}) and  the summation finally gets  easy to compute. We introduce here a notation that we shall use several times : 
\beqy \label{notationpairewisedistinct}
\Idis_n \ := \ \big\{ (i_1,\ld,i_n)\in \{1,\ld,N\}^n \ste  i_1,�\ld,i_n \trm{ are pairwise distinct}\big\}
\eeqy 
 (this set implicitly depends on $N$).

\subsection{Proof of b)} In this part, as $n$ is fixed, we shall denote $\bN_n$ (resp. $\bN_{-n}$) by  $\bMM=[M_{ij}]$ (resp. $\bK=[K_{ij}]$). For any non-negative integers $p,q,r,s$, wet set
\beq
m_{p,q,r,s} & := & \E \big( \tr  \bA^n \bMM \big)^p \ol{\big(\tr \bA^n \bMM \big)}^q \big(\tr  \bA^{-n} \bK  \big)^r \ol{\big( \tr  \bA^{-n} \bK \big)}^s   
\eeq
and our goal is to show that, as $N$ goes to infinity, the   numbers $m_{p,q,r,s}$ have limits satisfying 
conditions \eqref{00-112009}, \eqref{0-112009}, \eqref{3-112009}, \eqref{2-112009} and \eqref{4-112009} of the Lemma \ref{150112032015}. 
Note that  \eqref{00-112009}  and  \eqref{3-112009}  follow   from the fact the the Haar measure on the unitary group is invariant by multiplication by any $\mre^{\ii \tta} $, $\tta\in \R$. We shall use the following notations
$$
\begin{array}{rcll}
\ds\lim_{N \to \infty} \ff N \tr \big( \bMM \bMM^*\big)  =:   \al_\bMM &;& \
\ds\lim_{N \to \infty} \ds\ff N \tr \big( \bK \bK^*\big)  =: \al_\bK &; \
\vspace{2mm}\ds\lim_{N \to \infty} \ff N \tr \big( \bMM \bK\big)  =:  \bet_{\bMM,\bK} \\
\ds\lim_{N \to \infty} \ff N \tr \bMM =: \tau_\bMM &;& \ \ds\lim_{N \to \infty} \ff N \tr \bK =: \tau_\bK.
\end{array}
$$

\subsubsection{$\tr \bA^n \bMM $ and $\tr \bA^{-n}\bK$ are asymptotically two circular Gaussian complex variables satisfying conditions  \eqref{0-112009}  and \eqref{4-112009}}\la{section1022} We simply have to show that for any integer $p \geq 1$ 
\beqy\la{124151}
\E \big|\tr \bA^n \bMM \big|^{2p}                    & = & p! \big(b^{2n}((n-1)|\tau_\bMM|^2+\al_\bMM)\big)^p +\oo1\,, \\ \la{124152}
\E \big|\tr \bA^{-n}\bK \big|^{2p}                 & = & p! \big(a^{-2n}((n-1)|\tau_\bK|^2+\al_\bK)\big)^p +\oo1 \,,\\ \la{124153}
\E \big(\tr \bA^n \bMM  \ \tr \bA^{-n}\bK \big)^p  & = & p! \big((n-1) \tau_\bMM \tau_\bK + \bet_{\bMM,\bK}  \big)^p + \oo1\,.  
\eeqy
We shall prove it by induction on $p$. So first, we show the previous relation for $p=1$. Recall that we assume that $\bMM = \bV \wt\bMM \bV^*$ and $\bK = \bV \wt\bK \bV^*$, where $\wt\bMM$ and $\wt\bK$ are deterministic, so that, using the Lemma \ref{lem161423031500}, we have (denoting again by $\E_\bV$ the expectation with respect to the randomness of $\bV$),
\beq                                                     
\hspace{-15mm}\E_\bV{\big|\tr \bA^n \bV\wt\bMM\bV^* \big|^{2}}                           & = &  \ff N \tr \bA^n (\bA^*)^n  \left(\ff N \tr  \wt\bMM\wt\bMM^*  - \big|\ff N \tr \wt\bMM  \big|^2  \right)\\
&+ &\big| \tr \bA^n \big|^2\big|\ff N \tr \wt\bMM  \big|^2 + \OO{\ff N}\\
\hspace{-15mm}\E_\bV \big|\tr \bA^{-n}\bV\wt\bK\bV^* \big|^{2}                        & = &  \ff N \tr \bA^{-n} (\bA^*)^{-n} \left(\ff N \tr\big( \wt\bK\wt\bK^*\big) - \big|\ff N \tr \wt\bK  \big|^2  \right)\\
&+ &\big| \tr \bA^{-n} \big|^2\big|\ff N \tr \wt\bK  \big|^2 + \OO{\ff N}\\
\hspace{-15mm}\E_\bV \tr \bA^n \bV\wt\bMM\bV  \, \tr \bA^{-n}\bV\wt\bK\bV^*    & = & \ff N \tr \wt\bMM\wt\bK -\ff N \tr \wt\bMM \ff N \tr  \wt\bK  \\
&+&  \tr \bA^n \tr \bA^{-n} \ff N \tr \wt\bMM \ff{N}\tr \wt\bK + \OO{\ff N}\,.
\eeq        
This is asymptotically determined by the limits of  $\E{\big|\tr \bA^n  \big|^2}$, $\E{\big|\tr \bA^{-n} \big|^2}$, $\E \tr \bA^n  \tr \bA^{-n}  $,\\ $N^{-1}\E\tr \bA^{n} (\bA^*)^{n} $ and $N^{-1} \E\tr \bA^{-n} (\bA^*)^{-n} $. First, we compute  $\E\big|\tr  \bA^n  \big|^2$ for $n \geq 1$. We write
\beqy \la{12120303}
\E\big|\tr \bA^n  \big|^2 & = & \sum_{\ds^{1 \leq i_1,\ldots,i_n \leq N}_{1 \leq j_1,\ldots,j_n \leq N}} \Ec{u_{i_1 i_2} \cdots u_{i_n i_1} \ol{u_{j_1 j_2}} \cdots \ol{u_{j_n j_1}}} s_{i_1} s_{j_1} \cdots s_{i_n} s_{j_n}.
\eeqy
From \eqref{wg1}, we have a condition on the $i_k$'s and the $j_k$'s for a non-vanishing expectation, which is the multiset\footnote{We use the index $m$ in $\mset{\,\cdot\,}$ to denote a \emph{multiset}, which means that $\mset{x_1,\ldots,x_n}$ is the class of the $n$-tuple $(x_1,\ldots,x_n)$ under the action of the symmetric group $S_n$.} equality
\beqy \la{12000303}
\mset{i_1,\ldots,i_n} \ = \ \mset{j_1,\ldots,j_n},
\eeqy
The first consequence of \eqref{12000303} is that the sum is in fact over  $\OO{N^n}$   terms which all are at most   $\OO{N^{-n}}$, which means that any sub-summation over  $\oo{N^n}$   terms might be neglected. So from now on, we shall only sum over the  $n$-tuples  $(i_1,\ldots,i_n)\in \Idis_n$ (recall notation  \eqref{notationpairewisedistinct}).   Then \eqref{12120303} becomes
\beq
\E\big|\tr \bA^n  \big|^2 & = & \sum_{(i_1,\ld,i_n)\in\Idis_n} s_{i_1}^2\cdots s_{i_n}^2 \sum_{\sigma \in S_n} \Ec{u_{i_1 i_2} \cdots u_{i_n i_1} \ol{u_{i_{\sigma(1)} i_{\sigma(2)}}} \cdots \ol{u_{i_{\sigma(n)} i_{\sigma(1)}}}} + \oo1
\eeq
Let $c \in S_n$ be the cycle $(1 \; 2 \; \cdots n)$. From \eqref{wg1} (see Remark \ref{alldistinct}), as long as the $i_k$'s are pairwise distinct, one can write 
$$
\Ec{u_{i_1 i_2} \cdots u_{i_n i_1} \ol{u_{i_{\sigma(1)} i_{\sigma(2)}}} \cdots \ol{u_{i_{\sigma(n)} i_{\sigma(1)}}}} \ = \ \wg\big( \sigma c^{-1} \sigma^{-1} c \big)
$$
and from \eqref{wg2} and   Remark \ref{moebius}, we know that the non-negligible terms are the ones such that 
$
\sigma c^{-1} \sigma^{-1} c \ = \ id$, \ie $  \sigma c \ = \ c \sigma,
$
which means that $\sigma$ must be a power of $c$ and so, only $n$ permutations $\sigma$ contribute to the non negligible terms. At last, as $\Moeb(id) = 1$, we have
\beq
\E\big|\tr  \bA^n  \big|^2 & = & \sum_{(i_1,\ld,i_n)\in \Idis_n} s_{i_1}^2\cdots s_{i_n}^2 \ti n N^{-n}\big(1 + \oo1 \big)+ \oo1 \\
& = & n \left( \ff N \sum_{i=1}^N s_i^2 \right)^n + \oo1 \ = \ n b^{2n} + \oo1.
\eeq
In the same way, one can get
\beq
\E\big|\tr \bA^{-n} \big|^2 & = & n a^{-2n} + \oo1, \\
\E \tr \bA ^n \tr \bA^{-n}   & = & n + \oo1 .
\eeq
Let us now consider $N^{-1}\E \tr \bA^n (\bA^*)^n $ for $n\geq 1$. We have
\beqy \la{121203031544}
\ff{N}\E \tr \bA^n (\bA^*)^n  & = & N^{-1}\!\!\!\!\!\!\!\!\!\!\sum_{\sub{1 \leq i_0,i_1,\ld,i_{n} \leq N \\ 1 \leq j_0,j_1,\ld,j_{n} \leq N \\ i_0 = j_0, \ i_{n} = j_{n}}} \!\!\!\!\!\!\!\!\!\!\Ec{u_{i_0 i_1} \cdots u_{i_{n-1} i_{n}} \ol{u_{j_{0} j_{1}}} \cdots \ol{u_{j_{n-1} j_{n}}}} s_{i_1}s_{j_1}\cdots s_{i_{n}}s_{j_{n}}
\eeqy
As previously, we know that by \eqref{wg1}, that  for the    expectation to be non zero, we must have the multiset equality
\beqy \la{120003031544}
\mset{i_0,\ldots,i_n} \ = \ \mset{j_0,\ldots,j_n},
\eeqy
The first consequence of \eqref{120003031544} is that the sum is in fact over  $\OO{N^{n+1}}$   terms which  are all at most   $\OO{N^{-n-1}}$, so   that any sub-summation over  $\oo{N^{n+1}}$   terms might be neglected. As previously, we shall sum over the pairwise distinct indices $\Idis_{n+1}$ (see notation \eqref{notationpairewisedistinct}). Hence \eqref{121203031544} becomes
\beq
 N^{-1} \E \tr \bA^n (\bA^*)^n   & = & N^{-1} \!\!\!\!\!\!\!\!\!\!\sum_{(i_0,i_1,\ld,i_n)\in\Idis_{n+1}}\!\!\!\!\!\!\!\!\!\!s_{i_1}^2\cdots s_{i_{n}}^2 \sum_{\sub{\sigma \in S_{n+1} \\ \sigma(0)=0 \\ \sigma(n)=n}} \Ec{u_{i_0 i_1} \cdots u_{i_{n-1} i_{n}} \ol{u_{i_{\sigma(0)} i_{\sigma(1)}}} \cdots \ol{u_{i_{\sigma(n-1)} i_{\sigma(n)}}}} 
\eeq
Let $c \in S_{n+1}$ be the cycle $(0\; 1 \; 2 \; \cdots n)$. From \eqref{wg1} (see Remark \ref{alldistinct}) one can write 
$$
\Ec{u_{i_0 i_1} \cdots u_{i_{n-1} i_{n}} \ol{u_{i_{\sigma(0)} i_{\sigma(1)}}} \cdots \ol{u_{i_{\sigma(n-1)} i_{\sigma(n)}}}} \ = \ \wg\big( \sigma c^{-1} \sigma^{-1} c \big).
$$
As previously,  $\sigma$ must be a power of $c$ for the asymptotic contribution to be non-negligible. However, this time, we impose   $\sigma(0)=0$ and $\sigma(n)=n$, so that the only possible choice is $\sigma = id$ which means that only one term   contributes this time. At last,
\beq
\ff N \E \tr \bA^n (\bA^*)^n  & = & b^{2n} + \oo1,
\eeq
The same way, one can get
\beq
\ff N \E \tr \bA^{-n} (\bA^*)^{-n}  & = & a^{-2n} + \oo1.
\eeq
This concludes the first step of the induction.

In the second step, we have to prove the following induction relation: for any $p \geq 2$,
\beqy\la{313151}                                                     
\E \big|\tr \bA^n \bMM \big|^{2p}                  & = & p\E \big|\tr \bA^n \bMM \big|^{2}       \E \big|\tr \bA^n \bMM \big|^{2(p-1)}                        +\oo1 \\ \la{313152}
\E \big|\tr \bA^{-n}\bK \big|^{2p}                 & = & p\E \big|\tr \bA^{-n}\bK \big|^{2}    \E \big|\tr \bA^{-n}\bK \big|^{2(p-1)}           + \oo1 \\         \nonumber 
\E  (\tr \bA^n \bMM    \tr \bA^{-n}\bK  )^p   & = & p\Ec{\tr \bA^n \bMM  \ \tr \bA^{-n}\bK } \E\big[\big(\tr \bA^n \bMM  \ \tr \bA^{-n}\bK \big)^{p-1}\big]\\    &&\qquad \qquad\qquad\qquad\qquad\qquad\qquad\qquad\qquad\qquad + \oo1\la{313153} 
\eeqy                                                   
Let us   first  consider $\E \big|\tr \bA^n \bMM \big|^{2p}$. We shall use the following notation
\beq
\tr \bA^n \bMM & = & \sum_{i_0,i_1,\ldots,i_n} u_{i_0 i_1}s_{i_1} \cdots u_{i_{n-1} i_n} s_{i_n} M_{i_n i_0}  \ =: \ \sum_{\bi}  u_{\bi}s_{\bi}M_{i_{n} i_0}\, ,
\eeq
where the bold letter $\bi$ denotes the $(n+1)$-tuple $(i_0,\ld,i_n)$ and where we set \be\la{eq:notation31315}u_{\bi} := u_{i_0 i_1}\cdots u_{i_{n-1} i_n}\qquad ;\qquad s_{\bi} := s_{i_1}\cdots s_{i_n}\,.\ee 
Hence,
\beqy \la{11591200207003}
 \E \big|\tr \bA^n \bMM \big|^{2p}    \!\!\!\!  & = & \!\!\!\!\!\!\sum_{\substack{\bi^1,\ld,\bi^p\\ \bj^1,\ldots,\bj^p}}  \!\!\! \Ec{u_{\bi^1}\cdots u_{\bi^p} \ol{u_{\bj^1}}\cdots \ol{u_{\bj^p}}} s_{\bi^1}M_{i^1_n i^1_0}s_{\bj^1}\ol{M_{j^1_n j^1_0}}\cdots s_{\bi^p}M_{i^p_n i^p_0}s_{\bj^p}\ol{M_{j^p_n j^p_0}}
\eeqy
As usual, we know we can sum over the $\bi^k$'s satisfying that $(\bi^1,\ld,\bi^p)$ (the $p(n+1)$-tuple obtained by concatenation of the $\bi$'s) has  pairwise distinct entries and \st we have the set equality:
\beqy
  \rset{i^{\lam}_{\mu}, \ 1 \leq \lam \leq p,\ 0 \leq \mu \leq n}
 & =& \  \rset{j^{\lam}_{\mu}, \ 1 \leq \lam \leq p,\ 0 \leq \mu \leq n}.
\eeqy
Then, in order to have $\wg(id)$, we must match each of the  $(n+1)$-tuples $\bi^1,\ld,\bi^p$ with one of the  $(n+1)$-tuples  $\bj^1,\ldots,\bj^p$, \ie that for all $1 \leq \lam \leq p$, there is a $1 \leq \lam' \leq p$ such that we have the set  equality
$$
\{\bi^\lam\} \ =: \ \rset{i^{\lam}_{0},i^{\lam}_{1},\ld,i^{\lam}_{n}} \ = \ \rset{j^{\lam'}_0,j^{\lam'}_1,\ld,j^{\lam'}_n} \ := \ \{\bj^{\lam'}\}. 
$$
We rewrite \eqref{11591200207003} by summing according the possible choice to match $\{\bi^1\} = \rset{i^{1}_0,i^{1}_1,\ld,i^{1}_n}$
\beq
 \E \big|\tr \bA^n \bMM \big|^{2p}     & = & \sum_{\lam=1}^p  \sum_{\sub{\ds^{(\bi^1,\ld,\bi^p)\in \Idis_{p(n+1)}}_{(\bj^1,\ldots,\bj^p)\in\Idis_{p(n+1)}} \\ \bi^1 \leftrightarrow \bj^{\lam}}} \!\!\!\!\!\!\!\!\! \Ec{u_{\bi^1}\cdots u_{\bi^p} \ol{u_{\bj^1}}\cdots \ol{u_{\bj^p}}} s_{\bi^1}M_{i^1_n i^1_0}s_{\bj^1}\ol{M_{j^1_n j^1_0}}\cdots s_{\bi^p}M_{i^p_n i^p_0}s_{\bj^p}\ol{M_{j^p_n j^p_0}}+\oo1,
\eeq
where $\bi^1 \leftrightarrow \bj^\lam$ stands for the set equality 
$
\ds \rset{i^{1}_0,i^{1}_1,\ld,i^{1}_n}=\rset{j^{\lam}_0,j^{\lam}_1,\ld,j^{\lam}_n}
$. Then, we know that the set of indices $\rset{i^{1}_0,i^{1}_1,\ld,i^{1}_n}$ is disjoint from the others, so that by  Corollary \ref{lemsauveur15460502},
$$
\Ec{u_{\bi^1}\cdots u_{\bi^p} \ol{u_{\bj^1}}\cdots \ol{u_{\bj^p}}} \ = \ \Ec{u_{\bi^1}\ol{u_{\bj^{\lam}}}}\Ec{{u_{\bi^2}\cdots u_{\bi^p} \ol{u_{\bj^1}}\cdots\ol{u_{\bj^{\lam-1}}}\ol{u_{\bj^{\lam+1}}}\cdots \ol{u_{\bj^p}}}}
$$
and up to a proper relabeling of the indices, all the choices lead to the same value of the expectation, so that
\beq
&& \E \big|\tr \bA^n \bMM \big|^{2p}    \\
&&\hspace{-15mm} \ = \ p\sum_{\sub{\bi^1 \in \Idis_{n+1} \\ \bj^1 \in \Idis_{n+1}}}\!\!\! \Ec{u_{\bi^1}\ol{u_{\bj^1}}}s_{\bi^1}M_{i^1_n,i^1_0}s_{\bj^1}\ol{M_{j^1_n,j^1_0}}\!\!\!\!\!\!\!\!\!\sum_{\sub{(\bi^2,\ld,\bi^p)\in\Idis_{(p-1)(n+1)} \\ (\bj^2,\ldots,\bj^p)\in\Idis_{(p-1)(n+1)}}}\!\!\!\!\!\!\!\!\! \!\!\!\!\!\!\Ec{u_{\bi^2}\cdots \ol{u_{\bj^p}}} s_{\bi^2}M_{i^2_n i^2_0}s_{\bj^2}\ol{M_{j^2_n j^2_0}}\cdots s_{\bi^p}M_{i^p_n i^p_0}s_{\bj^p}\ol{M_{j^p_n j^p_0}}+\oo1 \\
&&\hspace{-15mm} \ = \ p\Ec{\big|\tr(\bA^n \bMM)\big|^{2}} \Ec{\big|\tr(\bA^n \bMM)\big|^{2(p-1)}} +\oo1.
\eeq
This proves \eqre{313151}. In the same way, we prove \eqre{313152} and \eqre{313153}, and thus conclude the proof of the induction.   \\
 
\begin{rmk}
In the last computation, we split the expectation and so we separated the summation implying that
$$
\Idis_{p(n+1)} \ = \ \Idis_{n+1} \ti \Idis_{(p-1)(n+1)}
$$
which is obviously inaccurate. Nevertheless, we easily see that
$$
\card \Idis_{p(n+1)}   \ = \ \card \lf(\Idis_{n+1} \ti \Idis_{(p-1)(n+1)}\ri) \big( 1 + \oo1\big), 
$$ 
which means that this inaccuracy is actually contained in the $\oo1$.
\end{rmk}

To conclude the proof of {\bf b)}, we have to prove that $\tr \bA^n \bMM $ and $\tr \bA^{-n}\bK $ satisfy Condition  \eqref{2-112009}  at the large $N$ limit.

\subsubsection{$\tr \bA^n \bMM $ and $\tr \bA^{-n}\bK $ satisfy Condition  \eqref{2-112009} at the large $N$ limit}
We apply the same idea as previously, but for a slightly   more complicated expectation. Let $p,q,r,s$ be positive integers and such that $p-q=r-s$. We denote joint moments by $m_{p,q,r,s}$:
\beqy\label{1211020420158741}
m_{p,q,r,s} & := & \E \big( \tr\bA^n \bMM \big)^p \ol{\big(\tr \bA^n \bMM \big)}^q \big(\tr \bA^{-n}\bK  \big)^r \ol{\big( \tr \bA^{-n}\bK \big)}^s   ,
\eeqy
and as
\beq
\tr \bA^n \bMM  & = & \sum_{i_0,i_1,\ldots,i_n} u_{i_0 i_1}s_{i_1} \cdots u_{i_{n-1} i_n} s_{i_n} M_{i_n i_0}  \ = \ \sum_{\bi}  u_{\bi}s_{\bi}M_{i_{n} i_0} \\
\tr \bA^{-n} \bK & = & \sum_{i_n,i_{n-1},\ldots,i_0} \ol{u}_{i_{n-1} i_{n}}s^{-1}_{i_n} \cdots \ol{u}_{i_{0} i_1} s_{i_1} K_{i_0 i_n}  \ = \ \sum_{\bi}  \ol{u}_{\bi}s^{-1}_{\bi}K_{i_{0} i_n} ,
\eeq
we rewrite \eqref{1211020420158741} as follow
\beqy\la{93100201503274bis}
\E  \sum_{\sub{\bi^1,\ld,\bi^{p} \\ \bj^1,\ld,\bj^{q} \\ \bk^1,\ld,\bk^{r} \\ \bl^1,\ld,\bl^{s} }} \prod_{\sub{1 \leq \lam \leq p \\ 1 \leq \mu \leq q \\ 1 \leq \nu \leq r \\ 1 \leq \tta \leq s}} \f{s_{\bi^{\lam} }s_{\bj^{\mu}}}{s_{\bk^{\nu}} s_{\bl^\tta}}  M_{i^\lam_n,i^\lam_0}  \ol{M}_{j^\mu_n,j^\mu_0} K_{k^\nu_0,k^\mu_n}\ol{K}_{\ell^\tta_0,\ell^\tta_n}  u_{\bi^{\lam}} u_{\bl^\tta} \ol{u_{\bj^{\mu}}} \ol{u_{\bk^{\nu}}}
\eeqy
(recall that the $s_{\bi}=s_{i_1}\cd s_{i_n}$ for $\bi=(i_0, \ld, i_n)$).
As previously, we deduce from Proposition \ref{wg} that for the non vanishing expectations, we must have the following multiset equality
\beqy\la{0204200151111}
&&\mset{i^{\lam}_\mu, \ 1 \leq \lam \leq p, 0 \leq \mu \leq n}\vcup\mset{\ell^{\lam}_\mu, \ 1 \leq \lam \leq s, 0 \leq \mu \leq n}   \nonumber \\  & =& \ \mset{j^{\lam}_\mu, \ 1 \leq \lam \leq q, 0 \leq \mu \leq n}\vcup\mset{k^{\lam}_\mu, \ 1 \leq \lam \leq r, 0 \leq \mu \leq n},
\eeqy
from which we  deduce that we can restrict the summation to the tuples   such that $$(\bi^1,\ld,\bi^p,\bl^1,\ld\bl^s)\in \Idis_{(p+s)(n+1)}$$  and that, for the non negligible terms (\ie those which lead to $\wg(id)$), we must match each of the $(n+1)$-tuples involved in $u$ (the $\bi$'s and the $\bl$'s) with one of those involved in $\ol u$ (the $\bj$'s and the $\bk$'s). For example, we sum according to the choice the ``partner'' of $\bi^1$.
\beq
 m_{p,q,r,s} 
&=& \sum_{\al=1}^q \E  \sum_{\sub{(\bi,\bl)\in\Idis \\ (\bj,\bk) \in \Idis \\ \bi^1 \leftrightarrow \bj^\al}} \prod_{\sub{1 \leq \lam \leq p \\ 1 \leq \mu \leq q \\ 1 \leq \nu \leq r \\ 1 \leq \tta \leq s}} \f{s_{\bi^{\lam} }s_{\bj^{\mu}}}{s_{\bk^{\nu}} s_{\bl^\tta}}  M_{i^\lam_n,i^\lam_0}  \ol{M}_{j^\mu_n,j^\mu_0} K_{k^\nu_0,k^\mu_n}\ol{K}_{\ell^\tta_0,\ell^\tta_n}  u_{\bi^{\lam}} u_{\bl^\tta} \ol{u_{\bj^{\mu}}} \ol{u_{\bk^{\nu}}} \\
&+&  \sum_{\bet=1}^r \E  \sum_{\sub{(\bi,\bl)\in\Idis \\ (\bj,\bk) \in \Idis \\ \bi^1 \leftrightarrow \bk^\bet}} \prod_{\ds^{\ds^{1 \leq \lam \leq p}_{1 \leq \mu \leq q}}_{\ds^{1 \leq \nu \leq r}_{1 \leq \tta \leq s}}} \f{s_{\bi^{\lam} }s_{\bj^{\mu}}}{s_{\bk^{\nu}} s_{\bl^\tta}}  M_{i^\lam_n,i^\lam_0}  \ol{M}_{j^\mu_n,j^\mu_0} K_{k^\nu_0,k^\mu_n}\ol{K}_{\ell^\tta_0,\ell^\tta_n}  u_{\bi^{\lam}} u_{\bl^\tta} \ol{u_{\bj^{\mu}}} \ol{u_{\bk^{\nu}}} + \oo1 .\\
\eeq
where, to simplify, $(\bi,\bl)$ stands for the $(p+s)(n+1)$-tuple obtained by the concatenation of the $\bi^\lam$'s and the $\bl^\mu$'s, and $\Idis$ implicitly means $\Idis_{(p+s)(n+1)}$.
As previously, we use the Corollary \ref{lemsauveur15460502} to split the expectations. Hence, one easily gets  
\beq
m_{p,q,r,s} & = & q \Ec{\big|\tr \bA^n \bMM \big|^{2}}m_{p-1,q-1,r,s} + r \Ec{\tr \bA^n \bMM \tr \bA^{-n}\bK } m_{p-1,q,r-1,s}+\oo1.
\eeq
To get the other relations, we just sum according to the choice of the partner of $\bj^1$ (resp. $\bk^1$ and $\bl^1$).

\subsection{Proof of a): asymptotic factorisation of joint moments}
The proof   relies mostly on   Corollary \ref{lemsauveur15460502}. We first   expand the expectation 
\beqy
\E\Big[\prod_{1 \leq n \leq n_0} \big(\tr \bA^{n}\bN_n \big)^{p_{n}}\ol{\big(\tr \bA^{n}\bN_n \big)^{q_{n}}} \big(\tr \bA^{-n}\bN_{-n} \big)^{r_{n}}\ol{\big(\tr \bA^{-n}\bN_{-n} \big)^{s_{n}}}\Big] .
\eeqy
Let $M^{(n)}_{ij}$ denote the $(i,j)$-th entry of $\bN_n$ and recall that for  $\bi = (i_0,\ld,i_n)$, we set \beq
u_{\bi} := u_{i_0 i_1}\cdots u_{i_{n-1} i_n}\qquad ;\qquad s_{\bi} := s_{i_1}\cdots s_{i_n}\,.\eeq 
We get
\beq
\tr \bA^n \bN_n  & = & \sum_{i_0,i_1,\ldots,i_n} u_{i_0 i_1}s_{i_1} \cdots u_{i_{n-1} i_n} s_{i_n} M_{i_n i_0}^{(n)}  \ = \ \sum_{\bi}  u_{\bi}s_{\bi}M_{i_{n} i_0}^{(n)} \\
\tr \bA^{-n} \bN_{-n}  & = & \sum_{i_n,i_{n-1},\ldots,i_0} \ol{u}_{i_{n_1} i_{n}}s^{-1}_{i_n} \cdots \ol{u}_{i_{0} i_1} s_{i_1} M^{(-n)}_{i_0 i_n}  \ = \ \sum_{\bi}  \ol{u}_{\bi}s^{-1}_{\bi}M_{i_{0} i_n}^{(-n)} ,
\eeq
so that 
\beqy\la{93100201503274}
\E \prod_{1 \leq n \leq n_0} \sum_{\substack{\bi^{n,1},\ld,\bi^{n,p_{n}} \\ \bj^{n,1},\ld,\bj^{n,q_{n}} \\ \bk^{n,1},\ld,\bk^{n,r_{n}} \\ \bl^{n,1},\ld,\bl^{n,s_{n}}}  } \prod_{\sub{1 \leq \lam \leq p_{n} \\ 1 \leq \mu \leq q_{n} \\ 1 \leq \nu \leq r_{n} \\ 1 \leq \tta \leq s_{n}}} \f{s_{\bi^{n,\lam} }s_{\bj^{n,\mu}}}{s_{\bk^{n,\nu}} s_{\bl^{n,\tta}}}  M^{(n)}_{i^{n,\lam}_n i^{n,\lam}_0}  \ol{M}^{(n)}_{j^{n,\mu}_n j^{n,\mu}_0} M^{(-n)}_{k^{n,\nu}_0 k^{n,\mu}_n}\ol{M}^{(-n)}_{\ell^{n,\tta}_0 \ell^{n,\tta}_n}  u_{\bi^{n,\lam}} u_{\bl^{n,\tta}} \ol{u_{\bj^{n,\mu}}} \ol{u_{\bk^{n,\nu}}}
\eeqy
where we use bold letters such as $\bi^{n,\lam}$ to denote  $(n+1)$-tuples  $ (i^{n,\lam}_0,i^{n,\lam}_1,\ldots,i^{n,\lam}_n)$. We can use the same ideas as in \cite[Lemma 5.8]{FloJean} to state that the non-negligible terms of the sum must satisfy that for all $n$, there are as much $(n+1)$-tuples involved in   $u$ as in $\ol{u}$, which means that
$$
 p_{n} + s_{n} \ = \  q_{n}+r_{n},
$$ 
and that we must have the multiset equalities 
\beqy\la{0094327035102}
&&\bigcup_{n=1}^{n_0}\mset{i^{n,\lam}_\mu, \ 1 \leq \lam \leq p_n, 0 \leq \mu \leq n}\vcup\mset{\ell^{n,\lam}_\mu, \ 1 \leq \lam \leq s_n, 0 \leq \mu \leq n}   \nonumber \\  & =& \ \bigcup_{n=1}^{n_0} \mset{j^{n,\lam}_\mu, \ 1 \leq \lam \leq q_n, 0 \leq \mu \leq n}\vcup\mset{k^{n,\lam}_\mu, \ 1 \leq \lam \leq r_n, 0 \leq \mu \leq n}.
\eeqy
   We deduce that there are a $\OO{N^{\sum_n n( p_n+s_n)}}$ non zero terms in \eqref{93100201503274} and we can easily show that any subsum over a $\oo{N^{\sum_n n( p_n+s_n)}}$ is negligible so that for now on we shall sum over the non pairwise indices. Then, we know that we can neglect any expectation $\E_\bU$ which won't lead to $\wg(id)$ (see \eqref{wg2}) so that \eqref{0094327035102} becomes
\beq
\forall 1 \leq n \leq n_0, & & \rset{i^{n,\lam}_\mu, \ 1 \leq \lam \leq p_n, 0 \leq \mu \leq n}\vcup\rset{\ell^{n,\lam}_\mu, \ 1 \leq \lam \leq s_n, 0 \leq \mu \leq n} \\
 & =& \  \rset{j^{n,\lam}_\mu, \ 1 \leq \lam \leq q_n, 0 \leq \mu \leq n}\vcup\rset{k^{n,\lam}_\mu, \ 1 \leq \lam \leq r_n, 0 \leq \mu \leq n}.
\eeq
It follows that the set of indices involved in the expansion of the $\tr \bA^n \bN_n $, $\tr \bA^{-n} \bN_{-n} $, $\ol{\tr \bA^n \bN_n }$, $\ol{\tr \bA^{-n} \bN_{-n} }$, is disjoint from  the set of indices involved in the expansion of the $\tr \bA^m \bN_m $, $\tr \bA^{-m} \bN_{-m} $, $\ol{\tr \bA^m \bN_m }$, $\ol{\tr \bA^{-m} \bN_{-m} }$, as long as $n\neq m$. Therefore,   Corollary \ref{lemsauveur15460502}  allows to conclude the proof of {\bf a)}.

%

\subsection{Proof of  c): computation of the limit covariance}

Let $f,g$ be polynomials in $z$ and $z^{-1}$ and let $\bM,\bN$ be $N\ti N$ deterministic matrices 
 such that, as $N\to\infty$, 
$$
\ff N \tr \bM \ \lto  \ \tau \quad ;\quad \ff N \tr \bN \ \lto  \ \tau'
\quad ;\quad \ff N \tr \bM\bN^* \ \lto  \ \al \quad ;\quad\ff N \tr \bM\bN \ \lto  \ \bet.
$$
We need to check that the limits of both sequences $$\E (\tr f(\bA)\bM  - a_0(f) \tr\bM) (\ovl{\tr g(\bA)\bN  - a_0(g) \tr\bN})$$ and $$\E (\tr f(\bA)\bM  - a_0(f) \tr\bM) (\tr g(\bA)\bN  - a_0(g) \tr\bN)$$ are the ones given in the statement of Theorem \re{theoremgen}. Note that it suffices to compute the limits for $f=g$ and $\bM=\bN$. Indeed, using   the classical polarization identities  for $\bM$ and $\bN$, first for general  polynomials $f,g$, we reduce the problem to the case $\bM=\bN$. Then, we use  polarization identities again to reduce the problem to $f=g$.

Also, recall that since $\bA \eloi e^{i\tta}\bA$ for any deterministic $\tta$, we know that for any positive  distinct integers $p,q$, we have 
\beq
\E \tr \bA^p\bM \tr\bA^{-q}\bM   \ = \ \E \tr \bA^p\bM \ol{\tr \bA^q \bM }  \ = \ 0.
\eeq
It follows, using \eqre{124151}, \eqre{124152} and \eqre{124153}, that 
\beq
\E \lf|\tr f(\bA)\bM  - a_0(f) \tr\bM  \ri|^2  & = & \sum_{\underset{\neq 0}{m,n \in \z}} a_m(f)\ol{a_n(f)}   \E \tr \bA^m \bM  \ol{\tr \bA^n \bM } \\
& = & \sum_{n \geq 1} \Big(|a_n(f)|^2 \E\big|\tr  \bA^n \bM \big|^2 + |a_{-n}(f)|^2 \E\big|\tr \bA^{-n} \bM \big|^2\Big)\\
& \tto & \sum_{n \geq 1} \big(|a_n(f)|^2 b^{2n}  + |a_{-n}(f)|^2 a^{-2n}\big)\big((n-1)|\tau|^2+\al\big),\\
\E \lf(\tr f(\bA)\bM  -  a_0(f)\tr\bM\ri)^2  & = & \sum_{\underset{\neq 0}{m,n \in \z}} a_m(f){a_n(f)}   \E \tr \bA^m \bM  \tr \bA^n \bM   \\
& = & \sum_{n \geq 1}2a_n(f)a_{-n}(f)\E \tr\bA^n \bM \tr \bA^{-n}\bM  \\
& \tto & 2 \sum_{n \geq 1} a_n(f) a_{-n}(f)\big( (n-1)\tau^2 + \bet \big),
\eeq
which concludes the proof.

\section{Proof of Corollary \re{carpol174151}} It is easy to see that, for any $z\notin S$, we have $$\log|\det(z-\bA)|=\beg{cases} \ds\Tr\log\bT+\real\mc{A}^N_z, & \trm{ with }\mc{A}^N_z:=\Tr\sum_{n\le -1}\f{\bA^{n}}{nz^n}\trm{ if $|z|<a$.}\\\ds N\log|z|+\real\mc{B}^N_z, & \trm{ with }\mc{B}^N_z:=-\Tr\sum_{n\ge 1}\f{\bA^n}{nz^n}\trm{ if $|z|>b$,}
 \en{cases}$$
 (in the first case, we used the fact that $|\det \bA|=\det\bT$). 
 Then, by Theorem \re{theoremgen}, $$\lf(\mc{A}^N_z\ri)_{|z|<a}\,\cup\, \lf(\mc{B}^N_z\ri)_{|z|>b}$$ converges, for the finite-dimensional   distributions, to  a centered complex Gaussian process  $$\lf(\mc{A}_z\ri)_{|z|<a}\,\cup\, \lf(\mc{B}_z\ri)_{|z|>b}$$ with covariance defined by $$\E \mc{A}_z\mc{A}_{z'}=0,\quad \E \mc{A}_z\ovl{\mc{A}_{z'}}=-\log(1-\f{z\ovl{z'}}{a^2}),$$
 $$\E \mc{B}_z\mc{B}_{z'}=0,\quad \E \mc{B}_z\ovl{\mc{B}_{z'}}=-\log(1-\f{b^2}{z\ovl{z'}}),$$
  $$\E \mc{A}_z\mc{B}_{z'} =-\log(1-\f{z'}{z}),\quad \E \mc{A}_z\ovl{\mc{B}_{z'}}=0,$$ where $\log$ denotes the canonical complex $\log$ on $B(1,1)$.
Then, one concludes by noting that for $A,B\in \C$, $2\real A\real B=\real(AB+A\ovl{B})$.

\appendix

\section{}
\subsection{Weingarten calculus on the unitary group} \la{sectionwg}
Here we summarize the   results we need about integration with respect to the Haar measure on unitary group, (see \cite[Cor. 2.4 and Cor. 2.7]{COL}). 

Let $\Moeb$ denote the   \emph{M\"obius function} of the lattice of non-crossing partition, defined for example 
in \cite[Lect. 10]{ns06}.
To simplify,  for any $k$-tuples $\bi =(i_1,\ldots,i_k)$ and $\bj=(j_1,\ldots,j_k)$, we set
$$
u_{\bi,\bj} \ := \ u_{i_1 j_1} u_{i_2 j_2}\cdots u_{i_k j_k}
$$
\begin{propo} \label{wg}
\indent Let $k$ be a positive integer and $\bU=(u_{ij})$ a $N \ti N$ Haar-distributed matrix. Let $\bi=(i_1,\ldots,i_k)$, $\bi'=(i'_1,\ldots,i'_k)$, $\bj=(j_1,\ldots,j_k)$ and $\bj'=(j'_1,\ldots,j'_k)$ be four $k$-tuples of $\left\{1,\ldots,N \right\}$. Then 
\begin{eqnarray} \label{wg1}
\E \left[ u_{\bi,\bj} \ol{u_{\bi',\bj'}}\right] \ = \ \sum_{\sigma, \tau \in S_k} \delta_{i_1,i'_{\sigma(1)}} \ldots \delta_{i_k,i'_{\sigma(k)}} \delta_{j_1,j'_{\tau(1)}} \ldots \delta_{j_k,j'_{\tau(k)}} \wg( \tau \sigma^{-1}),
\end{eqnarray}
where $\wg$ is a function called  the \emph{Weingarten function}. Moreover, for $\sigma \in S_k$,   the asymptotical behavior of $\wg(\sigma)$ is given by   
\begin{eqnarray} \label{wg2}
N^{k+|\sigma|}\wg(\sigma) &=& \Moeb(\sigma) + O\left( \frac{1}{N^2}\right),
\end{eqnarray}
where $|\sigma|$ denotes the minimal number of factors necessary to write $\sigma$ as a product of transpositions. 
\end{propo}
\begin{rmk} \label{alldistinct}
One should notice that if all $k$-tuples $(i_1,\ld,i_k)$, $(j_1,\ldots,j_k)$, $(i'_1,\ldots,i'_k)$,  and $(j'_1,\ldots,j'_k)$ have  pairwise distinct entries, then \eqref{wg1} becomes simpler because  in this case  there is at most one non-zero term in the sum.
\end{rmk}
\begin{rmk} \label{moebius}
\indent  The permutation $\sigma$ for which $\wg(\sigma)$ will have the largest order is the only one satisfying $|\sigma|=0$, i.e. $\sigma=id$. Also, $\Moeb(id) = 1$ (see \cite{COL}).
\end{rmk}
Here is a useful corollary which permits to simplify many computations.
\begin{cor} \label{lemsauveur15460502}
Let $\bi = (i_1,\ldots,i_p)$, $ \bj=(j_1,\ldots,j_p)$, $ \bk = (k_1,\ldots,k_q), \ \bl=(\ell_1,\ldots,\ell_q)$, $\bi' = (i'_1,\ldots,i'_p)$, $\bj'=(j'_1,\ldots,j'_p)$, $\bk' = (k'_1,\ldots,k'_q)$, $ \bl'=(\ell'_1,\ldots,\ell'_q) $ be tuples \st the multisets defined by $\bi$ and $\bi'$ (resp. by $\bj$ and $\bj'$, by $\bk$ and $\bk'$, by $\bl$ and $\bl'$) are   equal and  \st
$$
\{i_1, \ld, i_p\} \cap \{k_1,\ld, k_q\} \ = \ \{j_1, \ld, j_p\} \cap \{l_1,\ld, l_q\} \ =\ \emptyset  .
$$
Then
$$
\Ec{u_{\bi,\bj}u_{\bk,\bl}\ol{u_{\bi',\bj'}}\ol{u_{\bk',\bl'}}}  \ = \ \Ec{u_{\bi,\bj}\ol{u_{\bi',\bj'}}} \Ec{u_{\bk,\bl}\ol{u_{\bk',\bl'}}} \lf( 1 + \OO{\ff {N^2}} \ri).
$$
\end{cor}
\bpr
To prove this result, we first recall the exact expression of the M\"obius function : for any permutation $\sigma$ with cycle decomposition $C_1 C_2 \cdots C_r$, 
\beq
\Moeb(\sigma) & = & \prod_{i=1}^r (-1)^{|C_i|-1}\cat_{|C_i|-1},
\eeq
where $\cat_k$ is the $k$-th Catalan number, $ \ds\ff{k+1}\binom{2k}{k}$. Then, obviously, if $\sigma$ and $\tau$ are two permutations with disjoint supports, then
$$
\Moeb(\sigma\circ\tau) \ = \ \Moeb(\sigma)\Moeb(\tau) \quad \quad \text{ and } \quad \quad \big|\sigma\circ\tau\big| \ = \ |\sigma| + |\tau|.
$$
Thus
\beq
N^{p+q+|\sigma\circ\tau|}\wg(\sigma\circ\tau) & = & \Moeb(\sigma\circ\tau) + \OO{\ff{N^2}} \\
& = & \Moeb(\sigma)\Moeb(\tau) + \OO{\ff{N^2}} \\
& = & n^{p+|\sigma|} \wg(\sigma) N^{q + |\tau|} \wg(\tau) \lf(1 + \OO{\ff{N^2}}\ri)
\eeq
So that
$$
\wg(\sigma\circ\tau) \ = \ \wg(\sigma)  \wg(\tau) \lf(1 + \OO{\ff{N^2}}\ri).
$$
One can easily conclude.
\epr

We   also need the following lemmas in the paper.

\beg{lem}\la{Th1JAF}Let $\bA=\bU\bT$ with $\bU$ Haar-distributed on the unitary group and $\bT$ deterministic. 
Let $\eps>0$. There is a finite constant $C$ depending  only on $\eps$ (in particular, independent of $N$ and of $\bT$) \st for   all positive integer $n$ \st $n^6<(2-\eps)N$, we have
\be\la{59141}\E\Tr  \bA^{n}(\bA^{n})^*) \ \le \ CNn^2 \lf(m_2+\f{nm_\infty^2}{N}\ri)^{n}   \ee
and
\be\la{591412}\E[|\Tr  \bA^{n} |^2] \ \le \ C \lf(m_2+\f{nm_\infty^2}{N}\ri)^{n} ,  \ee
 where $m_2:=N^{-1}\Tr \bT\bT^*
$ and $m_\infty:=\|\bT\|_{\op{op}}$.
\end{lem}

\bpr See \cite[Th. 1]{BEN}.\epr

\begin{lem}\la{lem161423031500}
Let $\bV$ be Haar-distributed and let $\bA,\bB,\bC,\bD$ be deterministic $N\ti N$ matrices. Then we have 
\beq
\E \Tr  \bA\bV\bB\bV^*  \Tr \bC\bV\bD\bV^* & = & \ff{N^2-1}\lf(\Tr \bA\Tr\bC\Tr\bB\Tr\bD+\Tr \bA\bC\Tr\bB\bD \right) \\
& & \quad -\ff{N(N^2-1)}\lf(\Tr \bA\Tr\bC\Tr \bB \bD+\Tr \bA\bC\Tr \bB \Tr\bD\ri)
\eeq
\end{lem}

\bpr Let $ \mc{M}_N(\C)$ denote the set of $N\ti N$ complex matrices. It has already been proved, in \cite[Lem. 5.9]{FloJean}, that for any   matrices $\bA,\bB,\bC,\bD\in \mc{M}_N(\C)$, we have \beqy\label{LHT13813} \E\Tr \bA\bV\bB\bV^*\bC\bV\bD\bV^* &=&\ff{N^2-1}\lf\{\Tr \bA\bC\Tr\bB\Tr\bD+\Tr \bA\Tr\bC\Tr\bB\bD\ri\}\\ \nonumber&& -\ff{N(N^2-1)}\lf\{\Tr \bA\bC\Tr \bB \bD+\Tr \bA\Tr\bC\Tr \bB \Tr\bD\ri\}.\eeqy We deduce that  \beg{align}&\nonumber\E \bV\otimes\bV^*\otimes\bV\otimes\bV^*  = \ff{N^2-1}\lf\{\sum_{i,j,k,\ell} E_{j,k}\oti E_{k,j}\oti E_{i,\ell}\oti E_{\ell,i}+ E_{i,j}\oti E_{k,\ell}\oti E_{\ell,k }\oti E_{j,i} \ri\}\\ \label{LHT13813June2} &  \qquad   -\ff{N(N^2-1)}\lf\{\sum_{i,j,k,\ell} E_{j,k}\oti E_{\ell ,j}\oti E_{i,\ell}\oti E_{k,i}+ E_{i,j}\oti E_{j,k}\oti E_{k,\ell }\oti E_{\ell,i}\ri\},\end{align}
where the $E_{r,s}$ denote the elementary matrices. Indeed, the linear morphism $\Psi$ from   $ \mc{M}_N(\C)^{\otimes 4}$ to the space of $4$-linear forms on $ \mc{M}_N(\C)$ defined by $$\Psi(\bM\oti \bN\oti \mathbf{P}\oti \mathbf{Q})( \bA,\bB,\bC,\bD):=\Tr  \bA\bM\bB\bN\bC\mathbf{P}\bD \mathbf{Q}$$ is an isomorphism, and \eqre{LHT13813} proves that the left and right hand terms of \eqre{LHT13813June2} have the same image by $\Psi$. Then, applying  
$$\bM\oti \bN\oti \mathbf{P}\oti \mathbf{Q}\longmapsto   \Tr \big(\bA\bM\bB\bN\big)\Tr \big(\bC\mathbf{P}\bD\mathbf{Q}\big)$$ on both sides of \eqre{LHT13813June2}, we deduce the lemma.\epr

\subsection{Moments of a   Gaussian vector with values in $\C^2$} \label{11221031215}

The following lemma allows to prove that a complex random vector $(\tilde{X},\tilde{Y})$ is Gaussian without having to compute all its joint moments, by only proving an induction relation. 

\begin{lemme}\la{150112032015}
Let $(X,Y)$ be a  Gaussian random vector with values in $\C^2$ whose distribution is characterized by  
\beqy \label{00-112009}
\E X  \ = \ \E X^2  \ = \ \E Y  \ = \ \E Y^2  \ = \ 0,
\eeqy
and  
\beqy \label{0-112009}
\E |X|^2  = \sigma_X \   ;   \ \E |Y|^2  = \sigma_Y \ ; \ \E XY  = \sigma_{XY} \  ; \ \E X\ol Y  = 0.
\eeqy
Then, the moments 
\beqy \label{1-112009}
m_{p,q,r,s} & := & \E X^{p} \ol{X}^{q} Y^{r} \ol{Y}^{s} 
\eeqy 
satisfy
\beqy \label{3-112009}p-q \neq r - s   & \implies &   m_{p,q,r,s} \ = \ 0  ,
\eeqy
\beqy \label{2-112009}
m_{p,q,r,s} & = & \left\{\begin{array}{l}
q\sigma_X m_{p-1,q-1,r,s} + r\sigma_{XY} m_{p-1,q,r-1,s} \\
p\sigma_X m_{p-1,q-1,r,s} + s\ol{\sigma_{XY}} m_{p,q-1,r,s-1}\\ 
s \sigma_Y m_{p,q,r-1,s-1} + p\sigma_{XY} m_{p-1,q,r-1,s} \\
r \sigma_Y m_{p,q,r-1,s-1} + q\ol{\si_{XY}} m_{p,q-1,r,s-1} \\
\end{array}\ri. 
\eeqy
and
\beqy \label{4-112009}
m_{p,0,p,0}\;=\;\E (XY)^p  & = & p! \sigma_{XY}^p.
\eeqy
Conversely, if $(\tilde{X},\tilde{Y})$ is a    random vector with values in $\C^2$ \st both $\tilde{X}$ and $\tilde{Y}$ are Gaussian and  have joint moments $\tilde{m}_{p,q,r,s}$ satisfying  \eqref{00-112009}, \eqref{0-112009}, \eqref{3-112009}, \eqref{2-112009} and \eqref{4-112009}, then $(\tilde{X},\tilde{Y})$ is  Gaussian.
\end{lemme}

\bpr First, one easily obtains \eqref{3-112009} by noticing that for any $\tta\in \R$, $(\mre^{\ii\tta}X,\mre^{-\ii\tta}Y)\eqlaw (X,Y)$.

To prove the remaining, we consider  $(X,Y)$ as a real $4$-tuple $(\Re(X),\Im(X),\Re(Y),\Im(Y))=:(x_1,x_2,x_3,x_4)$ with covariance matrix
$$
\Gamma \ := \ \ff 2 \bpm \sigma_X & 0 & \Re(\sigma_{XY}) & \Im(\sigma_{XY}) \\
                         0 & \sigma_X & \Im(\sigma_{XY}) & -\Re(\sigma_{XY})  \\
                         \Re(\sigma_{XY}) & \Im(\sigma_{XY}) &  \sigma_Y & 0  \\    
                         \Im(\sigma_{XY}) & -\Re(\sigma_{XY}) & 0 & \sigma_Y   \\ 
\epm .
$$
Its Fourier transform is given, for $\bt = (t_1,t_2,t_3,t_4)$, by 
\beq
\Phi(\bt) & := & \E \exp\big(\ii(t_1 x_1 + t_2 x_2 + t_3 x_3 + t_4 x_4) \big)  \\
& = & \exp\big\{\ds - \ff 4\big( \sigma_X(t_1^2 + t_2^2) + \sigma_Y (t_3^2+t_4^2) \big) - \ff 2 \Re(\sigma_{XY})\big(t_1 t_3 - t_2 t_4 \big) - \ff 2 \Im(\sigma_{XY})\big(t_2 t_3 + t_1 t_4 \big)\big\} 
\eeq
We define the differential operators 
\beq
\partial_X \ = \ \partial_1 + \ii \partial_2 & \quad & \partial_{\ol X} \ = \ \partial_1 - \ii \partial_2 \\
\partial_Y \ = \ \partial_3 + \ii \partial_4 & \quad & \partial_{\ol Y} \ = \ \partial_3 - \ii \partial_4 \\
\eeq
so that
\beqy
\E X^p \ol{X}^q Y^r \ol{Y}^s  & = & (-\ii)^{p+q+r+s}\parmo{p}{q}{r}{s} 
\eeqy
and 
\beq
\partial_X \Phi(\bt) & = & -\ff 2 \big(\sigma_X (t_1 + \ii t_2) + \sigma_{XY}(t_3 - \ii t_4)\big)\Phi(\bt) \\
\partial_{\ol X} \Phi(\bt) & = & -\ff 2 \big(\sigma_X (t_1 - \ii t_2) + \ol{\sigma_{XY}}(t_3 + \ii t_4)\big)\Phi(\bt) \\
\partial_Y \Phi(\bt) & = & -\ff 2 \big(\sigma_Y (t_3 + \ii t_4) + \sigma_{XY}(t_1 - \ii t_2)\big)\Phi(\bt) \\
\partial_{\ol Y} \Phi(\bt) & = & -\ff 2 \big(\sigma_Y (t_3 - \ii t_4) + \ol{ \sigma_{XY}}(t_1 + \ii t_2)\big)\Phi(\bt) \\
\eeq
We easily deduce that
\beq
\parmo{p}{q}{r}{s}& = &  -q \sigma_X \parmo{p-1}{q-1}{r}{s} - r \sigma_{XY}\parmo{p-1}{q}{r-1}{s} 
\eeq
so that
\beq
\E X^p \ol{X}^q Y^r \ol{Y}^s  & = & q \sigma_X \E X^{p-1} \ol{X}^{q-1} Y^r \ol{Y}^s  + r \sigma_{XY}\E X^{p-1} \ol{X}^q Y^{r-1} \ol{Y}^s .
\eeq
This proves \eqre{2-112009}. 
To prove the last point, we simply write  
$$
\E (XY)^p   =  (-1)^p \partial^p_X \partial^p_Y \Phi(\bt)\big|_{ \bt=0}
\quad;\quad
\partial_X^p  \Phi(\bt)  =  \Big(-\f{\sigma_X (t_1 + \ii t_2) + \sigma_{XY}(t_3 - \ii t_4)}{2}\Big)^p \Phi(\bt).
$$
Then, using  Leibniz formula
and noticing that for all $k \leq p$, $$\partial^k_Y\Big(-\f{\sigma_X (t_1 + \ii t_2) + \sigma_{XY}(t_3 - \ii t_4)}{2}\Big)^p\big|_{\bt = 0} = \begin{cases} 0 \quad \text{ if } k < p \\ p! (\sigma_{XY})^p \text{ if } k=p \end{cases},$$
one can easily conclude.

Conversely, let   $(\tilde{X},\tilde{Y})$ be  a    random vector with values in $\C^2$ \st both $\tilde{X}$ and $\tilde{Y}$ are Gaussian and  have joint moments $ \tilde{m}_{p,q,r,s}   $ satisfying  \eqref{00-112009}, \eqref{0-112009}, \eqref{3-112009}, \eqref{2-112009} and \eqref{4-112009}. Let $\N$ denote the set of non-negative integers and let us define the sets   
$$
  \mathcal{K} : = \left\{ (p,q,r,s) \in \N^4\ste  \  p-q \neq r-s \ri\},
$$
$$
\mathcal{H}_0 : = \bigg\{ (p,q,r,s) \in \N^4\ste (r+s+|p-q|)\,(p+q+|r-s|)\,(p+r+|q-s|)\,(q+s+|p-r|)=0\bigg\}
$$
and, for $k\ge 1$, 
\beq
\mathcal{H}_k  & : = &  \big\{ (p,q,r,s)\in \N^4\ste \ (p-1,q-1,r,s) \text{ or } (p-1,q,r-1,s) \text{ or } \\
&& \quad \quad \quad  \quad \quad \quad\quad \quad \quad (p,q-1,r,s-1)\text{ or }(p,q,r-1,s-1) \in \mathcal{H}_{k-1}\big\}
\eeq
Then by hypothesis,   the joint moments function $\tilde{m}_{p,q,r,s} $ coincides with  $ m_{p,q,r,s} $ on 
 $  \mathcal{K}\cup  \mathcal{H}_0 $. Besides, by \eqref{2-112009},  if 
  $\tilde{m}_{p,q,r,s} $ coincides with  $ m_{p,q,r,s} $ on $\mathcal{H}_{k-1}$, then  $\tilde{m}_{p,q,r,s} $ coincides with  $ m_{p,q,r,s} $ on $\mathcal{H}_{k}$.
As
$$
\N^4 \ = \  \mathcal{K} \cup \bigcup_{k \geq 0} \mathcal{H}_k,
$$
we deduce that  $\tilde{m}_{p,q,r,s} $ coincides with  $ m_{p,q,r,s} $ on $\N^4$, wich implies that $(\tilde{X},\tilde{Y})\eqlaw(X,Y)$.
\epr

\noindent{\bf Acknowledgments.} The authors wish to thank the  anonymous referee for his careful reading and his useful advices.

  \begin{thebibliography}{10}

  \bibitem{greg-ofer} 
G. Anderson, O. Zeitouni \emph{A {CLT} for a band matrix model},  Probab. Theory Rel. Fields,
     2005, {134}, 283--338
  
     \bibitem{agz} G.~Anderson, A.~Guionnet, O.~Zeitouni \emph{An Introduction to Random Matrices}. Cambridge studies in advanced mathematics, {118} (2009).
     
     \bibitem{diaconis2003} A. D'Aristotile, P. Diaconis,   C. Newman \emph{Brownian motion and the classical groups}. With  Probability, Statisitica and their applications: Papers in Honor of Rabii Bhattacharaya. Edited by K. Athreya \emph{et al.} 97--116. Beechwood, OH: Institute of Mathematical Statistics, 2003.  
     
     \bibitem{baiysilver} Z.D.~Bai, J.~Silverstein \emph{CLT for linear spectral statistics of large-dimensional sample covariance matrices.} Ann. Probab. {32}, 2004, 533--605.

 \bibitem{BAI2009EJP} Z.D.~Bai,  X. Wang, W. Zhou  \emph{CLT for linear spectral statistics of Wigner matrices.} Electron. J. Probab. {14} (2009), no. 83, 2391--2417.
 \bibitem{BaiYaoBernoulli2005} Z.D.~Bai, J. Yao  \emph{On the convergence of the spectral empirical process of Wigner matrices.} Bernoulli {11} (2005) 1059--1092.
 
 \bibitem{BAI} J. Baik, G. Ben Arous,   S. P\'ech\'e \emph{Phase transition of the largest eigenvalue for nonnull complex sample covariance matrices}.  Ann. Probab., 33(5):1643--1697, 2005. 
     
       \bibitem{BD13} A. Basak, A. Dembo \emph{Limiting spectral distribution of sums of unitary and orthogonal matrices}. Electron. Commun. Probab. 18 (2013), no. 69, 19 pp.
   
\bibitem{bell} S.R. Bell \emph{The Cauchy transform, potential theory, and conformal mapping}. (1992) Studies in advanced Mathematics. CRC Press, Boca Raton, FL.

\bibitem{BENCLT} F. Benaych-Georges \emph{Central limit theorems for the Brownian motion on large unitary groups},
Bull. Soc. Math. France,Vol. 139, no. 4 (2011) 593--610.  

\bibitem{BEN} F. Benaych-Georges \emph{Exponential bounds for the support convergence in the Single Ring Theorem},  J. Funct. Anal.  268 (2015), pp. 3492--3507.

\bibitem{FloGuiJean} F. Benaych-Georges, G. C\'ebron, J. Rochet \emph{Fluctuation of matrix entries and application to outliers of elliptic matrices}, 	arXiv:1602.02929. 

\bibitem{AFCTCL} F. Benaych-Georges, A. Guionnet, C. Male \emph{Central limit theorems for linear statistics of heavy tailed random matrices}. 
Comm. Math. Phys. Vol. 329 (2014), no. 2, 641--686. 

\bibitem{BEN1} F. Benaych-Georges, R.N. Rao \emph{The eigenvalues and eigenvectors of finite, low rank perturbations of large random matrices}, Adv. Math. (2011), Vol. 227, no. 1, 494--521. 

\bibitem{BEN4} F. Benaych-Georges, R.N. Rao \emph{The singular values and vectors of low rank perturbations of large rectangular random matrices},  
J. Multivariate Anal., Vol. 111 (2012), 120--135. 

 \bibitem{FloJean} F. Benaych-Georges, J. Rochet \emph{Outliers in the Single Ring Theorem},  Probab. Theory Rel. Fields, Vol. 165 (2016), no. 1, 313--363. 
 
 \bibitem{BORCAP1} C. Bordenave, M. Capitaine \emph{Outlier eigenvalues for deformed i.i.d. random matrices}.
Comm. Pure Appl. Math. 69 (2016), no. 11, 2131--2194. 
 
 \bibitem{BOR1} C. Bordenave, D. Chafa\"{\i} \emph{Around the circular law}, Probab. Surv. 9 (2012), 1--89.
 
 \bibitem{b1906} \'E. Borel  \emph{Sur les principes de la th\'eorie cin\'etique des gaz}. Annales de l'\'Ecole Normale 
Sup\'erieure 23 (1906), 9--32.

\bibitem{CDFF}M. Capitaine, C. Donati-Martin, D. F\'eral, M. F\'evrier \emph{Free convolution with a semi-circular distribution and eigenvalues of spiked deformations of Wigner matrices}, Electron. J. Prob. Vol. 16 (2011), 1750--1792.

\bibitem{CebronKemp} G. C\'ebron, T. Kemp \emph{Fluctuations of Brownian Motions on $\mathbb{GL}_N$}. arXiv:1409.5624.

\bibitem{chatterjee} S.~Chatterjee
     \emph{Fluctuations of eigenvalues and second order {P}oincar\'e
              inequalities},
  {Probab. Theory Related Fields},
 {143},
     {2009},{1--40}.

     \bibitem{meckessourav08}  S. Chatterjee, E.  Meckes  \emph{Multivariate normal approximation using exchangeable pairs}  ALEA 4 (2008). 

\bibitem{mingo-piotr-collins-speicher07} B. Collins, J.A. Mingo,  P. \'Sniady,  R. Speicher  \emph{Second order freeness and fluctuations of random matrices. III. Higher order freeness and free cumulants}.  Doc. Math.  12  (2007), 1--70.
\bibitem{COL} B. Collins, P. Sniady \emph{Integration with respect to the Haar measure on unitary, orthogonal and symplectic group}.
Comm. Math. Phys., 264(3):773--795, 2006.

\bibitem{collins-stolz08} B. Collins, M. Stolz \emph{Borel theorems for random matrices from the classical compact symmetric spaces.} Ann.   Probab. Vol. 36, no. 3 (2008), 876--895. 

\bibitem{dia-shah} P. Diaconis, M. Shahshahani \emph{On the eigenvalues of random matrices},
Studies in applied probability, J. Appl. Probab. 31A (1994), 49--62.
\bibitem{dia-evans} P.~ Diaconis, S. Evans
     \emph{Linear functionals of eigenvalues of random matrices},
 {Trans. Amer. Math. Soc.},
{353}, {2001}, 2615--2633.

\bibitem{BRSVoverview}  B. Duplantier, R. Rhodes, S. Sheffield, V. Vargas \emph{Log-correlated Gaussian fields: an overview}, arXiv.

\bibitem{girko} V.L. Girko \emph{The elliptic law}, Teoriya Veroyatnostei i ee Primeneniya, 30(4):640--651, 1985.

    \bibitem{GUI} A. Guionnet, M. Krishnapur, O. Zeitouni \emph{The Single Ring Theorem}. Ann. of Math. (2) 174 (2011), no. 2, 1189--1217.
    
\bibitem{GUI2} A. Guionnet, O. Zeitouni \emph{Support convergence in the Single Ring Theorem}. Probab. Theory Related Fields 154 (2012), no. 3-4, 661--675.

  
 \bibitem{johansson} K.~Johansson \emph{On the fluctuations of eigenvalues of random Hermitian matrices.} Duke Math. J. {91} 1998, 151--204.
 
    \bibitem{jonsson}  D.~Jonsson \emph{Some limit theorems for the eigenvalues of a sample covariance matrix.} J. Mult. Anal. {12}, 1982, 1--38.
    
    \bibitem{jiang06}     T. Jiang  \emph{How many entries of a typical orthogonal matrix can be approximated by independent normals?} Ann. Probab. 34(4), 1497--1529. 2006.
    
    \bibitem{Kemp} T. Kemp \emph{Heat Kernel Empirical Laws on $\mathbb{U}_N$ and $\mathbb{GL}_N$}. J. Theor. Probab (2015). doi:10.1007/s10959-015-0643-7.
    
    \bibitem{KKP96}  A. M. Khorunzhy, B. A. Khoruzhenko, L. A. Pastur \emph{Asymptotic properties of large random
matrices with independent entries}, {J. Math. Phys.} 37 (1996) 5033--5060.  


\bibitem{thierrymylene} T. L\'evy, M. Ma\"{\i}da \emph{Central limit theorem for the heat kernel measure on the unitary group}. J. Funct. Anal. 259 (2010), no. 12, 3163--3204.

\bibitem{lytova} A.~Lytova,    L.~Pastur \emph{Central limit theorem for linear eigenvalue statistics of random matrices with independent entries}, {Ann. Probab.}, {37}, {2009}, 1778--1840.

\bibitem{meckes08}  E. Meckes  \emph{Linear functions on the classical matrix groups}.  Trans. Amer. Math. Soc.  360  (2008),  no. 10, 5355--5366.

\bibitem{mingo-nica04} J.A. Mingo, A. Nica   \emph{Annular noncrossing permutations and partitions, and second-order asymptotics for random matrices}.  Int. Math. Res. Not.  2004,  no. 28, 1413--1460. 

\bibitem{mingo-speicher06} J.A. Mingo, R. Speicher   \emph{Second order freeness and fluctuations of random matrices. I. Gaussian and Wishart matrices and cyclic Fock spaces}. J. Funct. Anal. 235 (2006), no. 1, 226--270.

\bibitem{mingo-piotr-speicher07}  J.A. Mingo, P.  \'Sniady, R. Speicher, \emph{Second order freeness and fluctuations of random matrices. II. Unitary random matrices}. Adv. Math. 209 (2007), no. 1, 212--240.
 \bibitem{ns06} A. Nica, R. Speicher \emph{Lectures on the combinatorics of free probability}. London Mathematical Society Lecture Note Series, 335. Cambridge University Press, Cambridge, 2006.

\bibitem{ROR12} S. O'Rourke, D. Renfrew \emph{Low rank perturbations of large elliptic random matrices}, 	Electron. J. Probab. 19 (2014), no. 43, 65 pp.

\bibitem{ORR} S. O'Rourke, D. Renfrew \emph{Central limit theorem for linear eigenvalue statistics of elliptic random matrices}, to appear in J. Theoret. Probab. 

\bibitem{ORRS} S. O'Rourke, D. Renfrew, A. Soshnikov \emph{On fluctuations of matrix entries of regular functions of Wigner matrices with non-identically distributed entries}. J. Theoret. Probab. 26 (2013), no. 3, 750--780.

\bibitem{SP06} S. P\'ech\'e \emph{The largest eigenvalue of small rank perturbations of Hermitian random matrices}, Prob. Theory Relat. Fields, 134 127--173, 2006.

\bibitem{peresvirag} Y. Peres, B.   Vir\'ag \emph{Zeros of the i.i.d. Gaussian power series: a conformally invariant determinantal process}. Acta. Math. 194. (2005) 1 -- 35.

\bibitem{PRS12} A. Pizzo, D. Renfrew, A. Soshnikov \emph{Fluctuations of matrix entries of regular functions of Wigner matrices}. J. Stat. Phys. 146 (2012), no. 3, 550--591. 

\bibitem{RiderJack} B. Rider, J. Silverstein \emph{Gaussian fluctuations for non-Hermitian random matrix ensembles}.
Ann. Probab. 34 (2006), no. 6, 2118--2143.

\bibitem{RiderVirag} B. Rider, B. Vir\'ag \emph{The noise in the circular law and the Gaussian free field}. Int. Math. Res. Not. IMRN 2007, no. 2.

     \bibitem{MShcherbina11} M.~Shcherbina
\emph{Central Limit Theorem for Linear Eigenvalue Statistics of the Wigner and
Sample Covariance Random Matrices}, Journal of Mathematical Physics, Analysis,
Geometry, 7(2), (2011), 176--192.

  \bibitem{sinai} Y.~Sinai,   A.~Soshnikov \emph{
     Central limit theorem for traces of large random symmetric
              matrices with independent matrix elements},
 {Bol. Soc. Brasil. Mat. (N.S.)}, {29}, {1998},
  1--24.

  \bibitem{soshni00} A. Soshnikov \emph{The central limit theorem for local linear statistics in
classical compact groups and related combinatorial identities}, Ann. Probab., 28 (2000), 1353--1370.

  \bibitem{TAO1} T. Tao 	\emph{Outliers in the spectrum of i.i.d. matrices with bounded rank perturbations}. Probab. Theory Related Fields 155 (2013), no. 1-2, 231--263.

 \en{thebibliography}

\end{document}